\documentclass[11pt,a4paper]{amsart}

\usepackage[bibtex]{ulthiel}

\title[The center of the asymptotic Hecke category]{The center of the asymptotic Hecke category and unipotent character sheaves}

\author[L. Rogel]{Liam Rogel}
\email{rogel@mathematik.uni-kl.de}

\author[U. Thiel]{Ulrich Thiel}
\email{ulrich.thiel@math.rptu.de}

\address{RPTU University Kaiserslautern-Landau, Department of Mathematics, 67653 Kaiserslautern, Germany}

\date{January 30, 2026}

\newaliascnt{construction}{theorem}
\newtheorem{construction}[construction]{Construction}
\aliascntresetthe{construction}

\providecommand{\abs}[1]{\lvert#1\rvert} 

\newcommand{\UCh}{\mathcal{U}}

\newcommand{\gbim}{\textnormal{-}\mathrm{gbim}}
\newcommand{\SBim}{\mathcal{H}}
\newcommand{\mbi}[1]{{\textbf{\em #1}}}

\newcommand{\PSL}{\mathrm{PSL}}

\usepackage{dsfont}
\usepackage{comment}

\usepackage{dynkin-diagrams}

\usepackage{todonotes}

\begin{document}

\begin{abstract}
 In 2015, Lusztig [Bull. Inst. Math. Acad. Sin. (N.S.)10(2015), no.1, 1–72] showed that for a connected reductive group over an algebraic closure of a finite field the associated (geometric) Hecke category admits a truncation in a two-sided Kazhdan--Lusztig cell, making it a categorification of the asymptotic algebra ($J$-ring), and that the categorical center of this ``asymptotic Hecke category'' is equivalent to the category of unipotent character sheaves supported in the cell. Subsequently, Lusztig noted that an asymptotic Hecke category can be constructed for any finite Coxeter group using Soergel bimodules. Lusztig conjectured that the centers of these categories are modular tensor categories (which was then proven by Elias and Williamson) and that for non-crystallographic finite Coxeter groups the $S$-matrices coincide with the Fourier matrices that were constructed in the 1990s by Lusztig, Malle, and Broué--Malle. If the conjecture is true, the centers may be considered as categories of ``unipotent character sheaves'' for non-crystallographic finite Coxeter groups.
 
 In this paper, we show that the conjecture is true for dihedral groups and for some (we cannot resolve all) cells of $H_3$ and $H_4$. The key ingredient is the method of $H$-reduction and the identification of the (reduced) asymptotic Hecke category with known categories whose center is already known as well. We conclude by studying the asymptotic Hecke category and its center for some infinite Coxeter groups with a finite cell.
\end{abstract}

\maketitle
\thispagestyle{empty}
\setcounter{tocdepth}{1}

\section{Introduction}

The representations of finite simple groups are a crucial ingredient in the investigation of finite symmetries. Most finite simple groups arise from finite reductive groups. An example of a reductive group is $G \coloneqq \mathrm{SL}_{n}(\overline{\mathbb{F}}_p)$ for a prime number $p$; its finite variants are $G(\mathbb{F}_{q}) \coloneqq \mathrm{SL}_{n}(\mathbb{F}_{q})$ for powers $q$ of $p$, yielding the finite simple groups $\PSL_n(\mathbb{F}_{q})$. An intrinsic construction produces from a reductive group $G$ a finite group $W$, the \emph{Weyl group}, which controls much of the structure of $G$. The Weyl group of $\mathrm{SL}_{n}(\overline{\mathbb{F}}_p)$, for example, is the symmetric group~$\mathfrak{S}_{n}$. Note that $W$ is independent of $p$.

Deligne--Lusztig theory \cite{DL-Reps} identifies an important subset of irreducible complex representations of finite reductive groups: the ``unipotent'' ones.
The work of Lusztig \cite{lusz1984charac} revealed an important feature of these representations: there is a finite set $U(W)$ just depending on $W$ which parametrizes the irreducible unipotent representations of $G(\mathbb{F}_q)$  \emph{independently} of $q$. Moreover, for each $\rho \in U(W)$ there is a polynomial $\mathrm{Deg}(\rho) \in \mathbb{Q} \lbrack \mbi{q} \rbrack$ such that when~$\mbi{q}$ is specialized to $q$, it gives the degree of the corresponding unipotent representation of $G(\mathbb{F}_q)$. Hence, the Weyl group controls the representation theory of the groups $G(\mathbb{F}_q)$.

Weyl groups admit a special presentation as abstract groups. For example, when taking the transpositions in $\mathfrak{S}_n$ as generators, the relations are the Artin braid relations. Coxeter groups are abstract groups admitting a more general such presentation. They still share many properties with Weyl groups. There is an especially well-behaved class of Coxeter groups: the crystallographic groups. Among the finite irreducible Coxeter groups, the crystallographic ones are precisely the Weyl groups; the non-crystallographic ones are the dihedral groups $I_2(n)$ (for $n=5$ and $n\geq 7$) and two further groups denoted $H_3$ and $H_4$.

The non-crystallographic groups cannot arise from a reductive group like Weyl groups do. Nonetheless, it has been observed that typical tools like the Hecke algebra, which is an algebra of equivariant functions on $G(\mathbb{F}_q)$ that again just depends on $W$, can be defined naturally for any Coxeter group \cite{Bourbaki-Lie}. In 1993 (\cite{lusztig1993} with an indication already in \cite{lusz1984charac}), Lusztig discovered the same phenomenon for unipotent representations: the sets $U(W)$ and the polynomials $\mathrm{Deg}(\rho)$ satisfy some natural properties; these still make sense when $W$ is (finite) non-crystallographic, and there are ad hoc constructions of such data satisfying these properties. Of course, there is no group $G(\mathbb{F}_q)$ where this could come from. Later, Lusztig \cite{lus94-exotic}, Malle \cite{malle94}, and Broué--Malle \cite{brouemalle93} found ad hoc constructions of Fourier matrices associated to (finite) non-crystallographic Coxeter groups. These are transition matrices between unipotent ``almost characters'' and unipotent characters for the groups $G(\mathbb{F}_q)$ introduced in \cite{lusz1984charac}. They are of fundamental importance since the almost characters have a uniform and simple construction. It is puzzling why there exist similar matrices when there is no reductive group. This observation was extended by Broué, Malle, and Michel \cite{BMM-Spetses} even to complex reflection groups~$W$. In the words of \cite{BMM-Spetses}, it almost looks like there is a ``fake algebraic group'' associated to complex reflection group. These mysterious objects were coined ``spetses''. Up to date, no one understands what they really are.

In 2015, Lusztig \cite{Lusztig2015} proposed a construction of categories associated with Coxeter groups. This construction conjecturally generalizes categories of unipotent character sheaves and naturally encodes the ad hoc Fourier matrices. If true, this would be an important step towards understanding spetses for non-crystallographic finite Coxeter groups. We will now summarize the basic line of thought towards the conjecture. 

\subsection{The conjecture and its background}

Let $G$ be a connected reductive group over $\overline{\mathbb{F}}_p$. Fix a prime $\ell \neq p$ and let $D^b_{G,c}(G)$ be the $G$-equivariant constructible bounded derived category of $\ell$-adic sheaves on $G$.
Lusztig's character sheaves \cite{CS1} are certain simple perverse sheaves in $D^b_{G,c}(G)$. 
As for characters, there is a notion of unipotent character sheaves \cite{CS1}. Taking the characteristic function of the Frobenius~$F$ on an unipotent character sheaf gives the corresponding unipotent almost character \cite{Shoji} for the finite group~$G^F$.\footnote{This only holds up to scalars and requires some restrictions. We can ignore this here.} The transition matrix $F_W$ between suitably normalized unipotent almost characters and unipotent characters is the Fourier matrix from \cite{lusz1984charac}. Like the parametrization of unipotent characters, it only depends on the Weyl group $W$ of $G$. 

Let $H_W$ be the Hecke algebra of $W$ with parameters as in \cite{EW-SoergelCalc}. The multiplicative properties of the Kazhdan--Lusztig basis $\lbrace b_w \rbrace_{w \in W}$ of $H_W$ lead to a decomposition of $W$ into two-sided cells \cite{KL-Hecke}.
Let $\UCh_G$ be the subcategory of $D^b_{G,c}(G)$ consisting of direct sums of unipotent character sheaves. To each unipotent character one can associate a unique two-sided cell $c$ of $W$ \cite{lusz1984charac} and this leads to a decomposition of $U(W)$ into subsets $U^c(W)$. This categorifies by \cite{CS3, Lusztig2015} and leads to a decomposition of $\UCh_G$ into subcategories $\UCh_G^c$. The Fourier matrix has block diagonal form with blocks $F_W^c$ indexed by the cells of $W$. 

Fix a Borel subgroup $B$ of $G$. Let $D^b_{B,c}(G/B)$ be the $B$-equivariant constructible bounded derived category of $\ell$-adic sheaves on $G/B$. This is a monoidal category with respect to convolution. The \emph{geometric Hecke category} $\mathcal{H}_G$ of $G$ is the subcategory of $D^b_{B,c}(G/B)$ consisting of semisimple perverse sheaves \cite{Williamson-HeckeCat}. This is a monoidal subcategory which categorifies $H_W$ in the sense that there is an algebra isomorphism
\begin{equation} \label{Thiel-Hecke-Categorification}
  H_W \to \lbrack \mathcal{H}_G \rbrack_\oplus \;, \quad b_s \mapsto \lbrack B_s \rbrack \;
\end{equation}
into the Grothendieck ring of $\mathcal{H}_G$, see \cite{RW-Tilting}. Here, $s \in W$ is a simple reflection and $B_s$ is the constant sheaf supported on $\overline{BsB}/B$. An important fact, which relies on the decomposition theorem \cite{BBD}, is that under this isomorphism the Kazhdan--Lusztig basis element $b_w$ for $w \in W$ gets mapped to a (uniquely characterized) direct summand $B_w$ of products of $B_s$ corresponding to a reduced expression of~$w$. This mirrors the properties of the Kazhdan--Lusztig basis and the indecomposable objects $\{ B_w \}_{w \in W}$ of $\mathcal{H}_G$ categorify the basis~$\{ b_w \}_{w \in W}$. 

Fix a cell $c$ and let $\mathcal{H}_G^c$ be the subcategory of $\mathcal{H}_G$ consisting of sheaves supported on $c$. This category is monoidal as well, but with respect to truncated convolution \cite{Lusztig-Cells}. We call it the \emph{asymptotic} Hecke category since it is a categorification of Lusztig's asymptotic algebra \cite{lusztig1987}. It follows from \cite{Lusztig-Cells, Bez2009, ost2014tensor} that there is a finite group $\Gamma_W^c$, a finite $\Gamma_W^c$-set $Y_W^c$, a 3-cocycle $\omega_W^c$ on $\Gamma_W^c$, and a monoidal equivalence
\begin{equation} \label{Thiel-Coh}
\mathcal{H}_G^c \simeq \mathrm{Coh}_{\Gamma_W^c}^{\omega_W^c}(Y_W^c \times Y_W^c) \;,
\end{equation}
the latter being the category of $\Gamma_W^c$-equivariant sheaves on $Y_W^c \times Y_W^c$ with convolution as tensor product and associator $\omega_W^c$. This description is key to understanding the following construction more explicitly. The \emph{Drinfeld center} of a monoidal category $\mathcal{C}$ is the category $\mathcal{Z}(\mathcal{C})$ of pairs $(Z,\gamma)$, where $Z \in \mathcal{C}$ and $\gamma$ is a functorial isomorphism
\begin{equation} \label{Thiel-Braiding}
  \gamma_X \colon X \otimes Z \overset{\simeq}{\longrightarrow} Z \otimes X
\end{equation}
for all $X \in \mathcal{C}$ which is compatible with the associator, see \cite[\S7.13]{EGNO}. Note that \eqref{Thiel-Braiding} induces a braiding on the Drinfeld center. From \eqref{Thiel-Coh} one obtains
\begin{equation} \label{Thiel-Vec}
\mathcal{Z}(\mathcal{H}_G^c) \simeq \Gamma_W^c\text{-}\mathrm{Vec}_{\Gamma_W^c}^{\omega_W^c}
\end{equation}
as braided monoidal categories, the latter being the category of $\Gamma_W^c$-equi\-variant $\Gamma_W^c$-graded vector spaces, see \cite[Example 8.5.4]{EGNO}. This is a modular tensor category \cite[\S8.13]{EGNO}, and it follows from \cite{lus94-exotic} that its $S$-matrix (which involves the braiding) is equal to the Fourier matrix $F_W^c$.

Lusztig \cite{Lusztig2015} gave geometric meaning to $\mathcal{Z}(\mathcal{H}_G^c)$ by constructing a mo\-no\-idal structure on $\UCh_G^c$ and establishing a natural monoidal equivalence
\begin{equation} \label{Thiel-center_equiv}
 \UCh_G^c \simeq \mathcal{Z}(\mathcal{H}_G^c) \;.
\end{equation}
In particular, $\UCh_G^c$ is a modular tensor category whose $S$-matrix is $F_W^c$. We note that the equivalence \ref{Thiel-center_equiv} seems to fit into a more general ``untruncated'' picture that is being established by the work of Bezrukav\-nikov--Finkel\-berg--Ost\-rik \cite{Bezrukavnikov2012}, Ben-Zvi--Nadler \cite{BenZvi2015}, and Bezrukavnikov\allowbreak--Ionov\allowbreak--Tolmachov\allowbreak--Varshavsky \cite{Bezrukavnikov2023}. \\

Let~$\mathfrak{h}$ be the root lattice of $G$. When placing $\mathfrak{h}$ in degree $2$ of the algebra $R$ of regular functions on $\overline{\mathbb{Q}}_\ell \otimes_{\mathbb{Z}} \mathfrak{h}$, then $R$ is as a graded algebra canonically isomorphic to $H_B^\bullet(\mathrm{pt},\overline{\mathbb{Q}}_\ell) \simeq R$, the total $B$-equivariant cohomology of $\overline{\mathbb{Q}}_\ell$ on a point. By \cite{RW-Tilting}, taking $B$-equivariant total cohomology on $G/B$ yields a fully-faithful monoidal graded functor
\begin{equation}
 \mathcal{H}_G \to R\gbim \;,
\end{equation}
the latter category being the category of graded $(R,R)$-bimodules. Hence, $\mathcal{H}_G$ is monoidally equivalent to a full subcategory of $R\gbim$: this is the category $\SBim_W$ of \emph{Soergel bimodules} introduced by Soergel \cite{Soergel-BimodulnOld, Soergel-Bimoduln}. The key feature of this category is that it can be constructed just from $W$ (and a reflection representation $\mathfrak{h}$). Moreover, it can be defined naturally for \emph{any} Coxeter group (when choosing an appropriate reflection representation $\mathfrak{h}$) and it yields a categorification of the Hecke algebra $H_W$ generalizing \eqref{Thiel-Hecke-Categorification}. It is a deep theorem by Elias and Williamson \cite{EW-Hodge} that the indecomposable objects $\{ B_w \}_{w \in W}$ of $\SBim_W$ categorify the Kazhdan--Lusztig basis as before. We should thus think of $\SBim_W$ as the ``Hecke category'' of spetses of type $W$. This category provides us with a kind of ``categorical geometry'' even if there is no reductive group.

For a two-sided cell $c$ in a finite Coxeter group $W$, Lusztig \cite[\S10]{Lusztig2015} defined an asymptotic Hecke category $\SBim_W^c$ and a monoidal structure on it, mimicking that of the asymptotic geometric Hecke category $\mathcal{H}_G^c$. Lusztig then took its Drinfeld center
\begin{equation}
 \UCh_W^c \coloneqq \mathcal{Z}(\SBim_W^c) \;.
\end{equation}
This should be considered as a category of ``unipotent character sheaves'' on the spetses of type $W$. Consequently, it should satisfy several properties. First of all, $\UCh_W^c$ should be a modular tensor category as conjectured by Lusztig \cite[\S10]{Lusztig2015}. This is indeed true and was proven by Elias--Williamson \cite{elias.will16}. We are thus down to the following conjecture.

\begin{conjecture}[Lusztig {\cite[\S10]{Lusztig2015}}] \label{lusztig_conjecture}
	Let $W$ be a non-crystallographic finite Coxeter group and let $c$ be a cell in $W$. The $S$-matrix of $\UCh_W^c$ is equal to the Fourier matrix $F_W^c$ from \cite{lus94-exotic, malle94}. In particular, the number of simple objects of $\UCh_W^c$ is equal to the number of unipotent characters supported in $c$.
\end{conjecture} 

The conjecture provides a \emph{uniform} categorification — and thus deeper meaning — of the ad hoc constructions of unipotent characters and the Fourier transform for non-crystallographic finite Coxeter groups. In fact, since unipotent characters are so important invariants of a group, maybe one can actually consider the categories $\UCh_W^c$ as ``truncated spetses'' for finite Coxeter groups. Following this idea, the actual spetses for finite Coxeter groups may be obtained as the center of the homotopy category of the Hecke category, which would fit in the work of Ben-Zvi--Nadler \cite{BenZvi2015}. Of course, the real spetses idea concerns complex reflection groups. A Hecke category for complex reflection groups still needs to be discovered, but once it is discovered one may proceed along the same path. 

Finally, the Fourier transform matrix for the big cell for $H_4$ given by Malle \cite{malle94} is not yet known to be an $S$-matrix of a modular tensor category: the conjecture provides, for the first time, a precise candidate. This may have implications in other fields (like physics) as well.

\subsection{Results in this paper}

First, we note that the asymptotic Hecke category $\SBim_W^c$ is \emph{multifusion} in the language of \cite{EGNO}, see \Autoref{sub:asymptotic_hecke}. It thus has a component fusion subcategory $\SBim_W^h$, corresponding to a diagonal $H$-cell~$h$ in $c$, see \Autoref{sub:h-reduction}. A crucial observation is that the centers of $\SBim_W^h$ and $\SBim_W^c$ are equivalent, see \Autoref{equ:h_reduction_center} and the general \Autoref{prop:centr_reduction}. We can thus work with $\SBim_W^h$, which is simpler. 

We show in \Autoref{sec:dihedral} (see \Autoref{thm:dihedral}) that \Autoref{lusztig_conjecture} holds for dihedral groups. The key is to identify $\SBim_W^h$ with the even part of the Verlinde category and noticing that the fusion data of the center of the latter categories are already in the literature. To be more precise, while the asymptotic Hecke algebra can be seen directly to be isomorphic to the Grothendieck ring of the even part of the Verlinde category, written $\mathrm{Ad}(\mathcal{C}_n)$, we need more known results to see that this algebra is not categorified by a different category, see \Autoref{sub:adjoint_category}. This allows us to compute with the asymptotic Hecke category without needing the category of Soergel bimodules. The center of $\mathrm{Ad}(\mathcal{C}_n)$ has also been computed in the literature but without any connections to our setting. We describe the computation, give small examples, and show how its $S$-matrix coincides with the Fourier matrix by Lusztig.

By similar means we confirm Conjecture \ref{lusztig_conjecture} for some (we cannot resolve all) cells of $H_3$ and $H_4$ in Section \ref{sec:h34}. In some cells we still have two different options for the categorification, and we point out which one is the ``right'' one assuming \Autoref{lusztig_conjecture} holds. Only the middle cell in type $H_4$, the one of $a$-value $6$, remains a complete mystery we cannot resolve yet. 

Finally, we note that the asymptotic Hecke category can also be constructed for arbitrary (not necessarily finite) Coxeter groups and a finite Kazhdan--Lusztig cell. In \Autoref{sec:infinite_cases} we study infinite Coxeter groups having a finite cell of $a$-value equal to or less than $2$. We describe the corresponding asymptotic Hecke category and its center. Even though we did not find new fusion or modular tensor categories in these examples, we expect some new examples will arise from the setting of asymptotic Hecke categories.

We begin in Section \ref{sec:asymptotic_hecke} with a detailed review of the construction of the asymptotic Hecke category. In Section \ref{sec:center_asymptotic} we discuss generalities about its center and summarize known results in the Weyl group case. For all except~$3$ so-called exceptional cells $c$ in type $E_7$ and $E_8$ the asymptotic Hecke category is known to be of the form $\mathrm{Coh}_{G_c}(X_c\times X_c)$ for some group $G_c$ and a $G_c$-set $X_c$. The possibilities for $(G_c,X_c)$ are due to Lusztig and listed in \Autoref{ex:classification}. We show in \Autoref{rem_center_hreduction} how the Drinfeld center of the multifusion category $\mathrm{Coh}_{G_c}(X_c\times X_c)$ is equivalent to that of $\mathrm{Vec}(G_c)$ using the method of $H$-reduction as described in \Autoref{sub:h-reduction}. The $S$-matrices are listed in \Autoref{cor_smatrix_weyl}, and we get the same matrices as the combinatorially computed results of Lusztig in \cite{lusz1984charac}.

\subsection*{Acknowledgements}
We would like to thank Ben Elias, Daniel Tubbenhauer, and Geordie Williamson for many helpful discussions on this topic. The first author further thanks Ben Elias for his hospitality during a two months stay at the University of Oregon last summer. We would furthermore like to thank Fabian Mäurer for the development of software for computing the center of a fusion category \cite{MaurerThiel-Center,Maurer-TensorCategories} which helped us find some key ideas. We would like to thank Gunter Malle for comments on a preliminary version of this paper. We also thank the anonymous referees of ``Representation Theory'' for their many helpful comments which greatly improved the paper. 
This work was supported by the SFB-TRR 195 ``Symbolic Tools in Mathematics and their Application'' of the German Research Foundation (DFG). 

\tableofcontents

\section{The asymptotic Hecke category} \label{sec:asymptotic_hecke}

We describe the construction of the asymptotic, or truncated, Hecke category categorifying the asymptotic Hecke algebra associated to a two-sided Kazhdan--Lusztig cell of a Coxeter group. The main construction is due to Lusztig \cite[Section 10]{Lusztig2015}.

We start with the construction of the asymptotic Hecke algebra and then go into detail into the construction of the asymptotic Hecke category.

\subsection{The asymptotic Hecke algebra} \label{sec:jring}

We use the notation of \cite{IntroSoergel}. For a Coxeter system $(W,S)$, where $S=\{s_i\}\subset W$ is the set of reflections generating $W$ via $W = \langle s_i \mid s_i^2=(s_is_j)^{m_{i,j}}=1, m_{i,j} \in \ZZ \cup \{\infty\} \rangle$, we denote by $H_W$ the \emph{(equal parameter) Hecke algebra}, a unital associative algebra over $A\coloneqq \ZZ[v^{\pm 1}]$ generated by elements $\delta_w$ for $w\in W$ and subject to the quadratic relation $(\delta_s-v^{-1})(\delta_s+v)=0$ and braid relation $\underbrace{\delta_s\delta_t\ldots}_{m_{s,t}} = \underbrace{\delta_t\delta_s\ldots}_{m_{t,s}}$ for $m_{s,t}<\infty$.  
The Kazhdan--Lusztig basis is denoted by $\{b_w\mid w\in W \}\subseteq H_W$, and we write 
\begin{equation}
	b_x = \delta_x + \sum_{y<x}{h_{y,x}\delta_y}
\end{equation}  
with the Kazhdan--Lusztig polynomials $h_{y,x}\in v\ZZ[v]$. Furthermore, we define polynomials $h_{x,y,z}$ in $A$ such that \begin{equation}
	b_xb_y=\sum_{z\in W}{h_{x,y,z}b_z}.
\end{equation} We write $z\leftarrow_L y$ if there exists an element $x\in W$ such that $h_{x,y,z}\neq0$. We extend this relation transitively to a preorder $\leq_L$, i.e. $z\leq_L y$ if there exists a sequence $z=x_0,x_1,\ldots,x_k=y$ with $x_i\leftarrow_L x_{i+1}$ for all $i$. Furthermore, $z\sim_L y$ if and only if both $z\leq_L y$ and $y\leq_L z$. We call an equivalence class with respect to $\sim_L$ a \emph{left} or \emph{$L$-(Kazhdan--Lusztig) cell}.

Similarly, we define the relation $z\leftarrow_R y$ if there is an $x\in W$ such that $h_{y,x,z}\neq0$, and extend this in the same way to an equivalence relation $\sim_R$ whose equivalence classes are called \emph{right} or \emph{$R$-cells}. Finally, let $x\sim_J y$ be the extension of the relation $x\leftarrow_J y$, which is defined as $x\leftarrow_L y$ or $x\leftarrow_R y$, and let the equivalence classes of $\sim_J$ be called $J$- or \emph{two-sided}-cells. These relations are due to Green \cite{Green1951OnTS} for monoids and have been extended to algebras and categories.
On $W$ we define the \emph{$a$-function} \begin{equation}
	a: W\to \NN\cup{\infty},
\end{equation} where for each $z\in W$, $a(z)$ is the smallest integer in $\NN$ such that \begin{equation}
	v^{a(z)}h_{x,y,z}\in \ZZ[v] \text{ for all } x,y\in W,
\end{equation} or $\infty$ if no such integer exists.
 It is conjectured that no case with an infinite value occurs, see \cite[Section 14.2 and 15]{Bonnafe18}. We will only consider \emph{bounded} Coxeter groups, i.e. of finite $a$-value. For any $x,y,z\in W$, we now define \begin{equation}
	\gamma_{x,y,z^{-1}}\coloneqq (h_{x,y,z}v^{a(z)})(0)\in \ZZ, 
\end{equation}to be the coefficient of the $v^{-a(z)}$-term in $h_{x,y,z}$. 

Using the coefficients $\gamma_{x,y,z}$ one defines a new ring structure on the set $\langle j_w \mid w \in W \rangle_{\ZZ}$, see \cite{lusztig1987}. The \emph{asymptotic Hecke algebra} or \emph{$J$-ring} $J\coloneqq J_W$ is the free abelian group generated by $\{j_x \mid x\in W \}$ subject to the relations \begin{equation}
	j_xj_y=\sum_{z\in W}{\gamma_{x,y,z^{-1}}j_{z}},  
\end{equation} for all $x,y,z$. 

\begin{definition}[{\cite[Section 3.1]{EGNO}}]\label{def:fusion_ring}
	Let $R$ be a unital ring which is free as a $\ZZ$-module. The tuple $(R,B)$ for a fixed basis $B=\{b_i\}_{i\in I}$ of $R$ is called a \emph{based ring} if we have: \begin{enumerate}
		\item $b_ib_j=\sum_{k\in I}{c_{i,j}^kb_k}$,  for $c_{i,j}^k\in \ZZ_{\geq0}$.
		\item The unit $1\in R$ is a non-negative linear combination of basis elements. Let $I_0\subset I$ index those $b_i$ occurring in the decomposition of $1$, and define $\tau:R\to\ZZ$ by $\tau(b_i)=1$ if $i\in I_0$ and $\tau(b_i)=0$ otherwise.
		\item There exists an involution $i\mapsto i^*$ on $I$ such that $b_i\mapsto b_{i^*}$ extends to a ring anti-involution on $R$, and $\tau(b_{i}b_j)=\delta_{j,i^*}$ (i.e., $b_ib_{i^*}$ contains exactly one basis summand of the unit, with coefficient 1).
	\end{enumerate}
	If the basis is finite, i.e. $R$ is of finite rank, we call it a \emph{multifusion ring}. If furthermore $1\in B$ we call it a \emph{fusion ring}.
\end{definition}

\begin{remark}\label{rem:duflo}
	For finite $W$ one always has $\gamma_{x,y,z}\geq 0$. This has been shown for crystallographic $W$ in \cite[Lemma 5.2(d)]{lusztig1985cells} and \cite[1.1(e)]{lusztig1987}, and verified by explicit calculation for non-crystallographic types $H_3$, $H_4$, and $I_2(p)$ in \cite{ducloux2006731,lusztig2014hecke}.
	The $J$-ring of a finite Coxeter group is a multifusion ring with unit $\sum_{t\in D}{j_t}$, where $D$ denotes the set of \emph{Duflo involutions} in $W$. For general bounded Coxeter groups, it is conjectured that the formal sum $\sum_{t\in D}j_t$ still acts as a local unit; see \cite[Section 13.4, Conjecture 14.2, Section 18.3]{lusztig2014hecke}.
\end{remark}

\begin{remark}\label{rem:directsum}
	By \cite[Corollary 1.9]{lusztig1987}, $\gamma_{x,y,z}\neq0$ implies $x,y,z$ lie in the same two-sided cell. Thus, defining $J_c\coloneqq \langle j_x\mid x\in c\rangle$ for each two-sided cell $c\subseteq W$, we obtain a decomposition
	\begin{equation}
		J_W \simeq \bigoplus_{c\subset W}J_c.
	\end{equation}
	We call $J_c$ the \emph{asymptotic Hecke algebra} associated to $c$. Each summand $J_c$ is a multifusion ring with unit given by the sum of all Duflo involutions in~$c$.
\end{remark}

\subsection{Construction of the asymptotic Hecke category} \label{sub:asymptotic_hecke}

Let now $\mathcal{H}_W$ be the category of Soergel bimodules associated to a given Hecke algebra $H_W$, see \cite{IntroSoergel}.

In \cite{elias.will16}, Elias and Williamson showed that the asymptotic Hecke category endowed with the monoidal product motivated by Lusztig in \cite[Section 10]{Lusztig2015} forms a rigid, see \cite[Section 2.10]{EGNO}, category. This implies that the asymptotic Hecke category for a two-sided cell with finitely many left cells is multifusion, see \Autoref{rem:multifusion_asymptotic}. For finite Weyl groups we list in \Autoref{ex:classification} for which cases a description of the asymptotic Hecke category is known. 

We go through their computations and motivate the construction of the asymp\-totic Hecke category in parallel to that of the asymptotic Hecke algebra. One key observation of Elias and Williamson is that the direct sum decomposition for Soergel bimodules is not canonical, therefore we get problems if we would naively try to define an asymptotic monoidal product by just taking the ordinary monoidal product and sending it to the lowest graded summand. Following \cite[Section 10]{Lusztig2015} one can define a canonical direct sum decomposition using the perverse filtration on Soergel bimodules.

\begin{example}
	This is seen in \cite[Example 2.1]{elias.will16}. For $s\in W$ a reflection and $B_s\in \mathcal{H}_W$ the Bott--Samelson bimodule corresponding to $s$, we have $B_s\otimes B_s\simeq B_s(+1)\oplus B_s(-1)$. In the Hecke algebra we have accordingly $b_sb_s=(v+v^{-1})b_s$. The $a$-value of $s$ is $1$ and in the $J$-ring this implies $j_sj_s=j_s$.
	
	One would like to find morphisms inside $\mathcal{H}_W$ from $B_s(-1)$ to $B_s\otimes B_s$ and vice versa to construct a categorification of the $J$-ring. However, by Soergel's Hom-formula, while the graded rank of the space $\mathrm{Hom}_{\mathcal{H}_W}(B_s\otimes B_s,B_s(-1))$ is $v^{-1}+2v+v^3$ and therefore a projection to $B_s(-1)$ is unique up to scalar, the inclusion is not unique. Two different direct sum decompositions can be found in \cite[Exercise 8.39 and 8.42]{IntroSoergel}.
	
	This means that $B_s(-1)$ is not a canonical subobject and one cannot directly replicate the multiplication of the $J$-ring on the category level.
\end{example} 

This shows that the lowest graded summand of a Soergel bimodule is not canonical. While the multiplication in the $J$-ring can be defined by ignoring all higher gradings, we cannot just define a monoidal product in the same way. The main result of  \cite{elias.will16} was to show relative hard Lefschetz for Soergel bimodules, as this allows to talk about certain ``canonical'' submodules of Soergel bimodules. To be more precise, we call a Soergel bimodule $B$ \emph{perverse} if it is isomorphic to a direct sum of Bott--Samelson bimodules without shifts. For an arbitrary Soergel bimodule $B$, the \emph{perverse filtration} is of the form \begin{equation}
	\ldots \subset \tau_{\leq i} B\subset \tau_{\leq i+1}B\subset \ldots,
\end{equation} where $\tau_{\leq i} B$ lies in the full subcategory of Soergel bimodules only generated by objects $B_x(m)$ for $m\geq -i$. Similarly, we consider $B/\tau_{\leq i} B$ which lies in the full subcategory of Soergel bimodules only generated by objects $B_x(m)$ for $m<-i$. We then write $H^i(B)\coloneqq (\tau_{\leq i}B/\tau_{< i}B)(i)$ for the \emph{perverse cohomology} of $B$.

\begin{theorem}[{\cite[Theorem 1.2]{elias.will16}}]
	Fix a Coxeter system $(W,S)$ and let $\mathcal{H}_W$ be the associated category of Soergel bimodules. If $\rho\in\mathfrak{h}^*$ is dominant regular (i.e. the Demazure operator is positive, $\partial_s(\rho)>0$, for all $s\in S$) and $x,y\in W$ are arbitrary the morphism \begin{equation}
		\eta: B_x\otimes_R B_y \to B_x\otimes_R B_y(2),~b\otimes b'\mapsto b\rho\otimes b'=b\otimes \rho b' 
	\end{equation} induces an isomorphism \begin{equation}
	\eta^i: H^{-i}(B_x\otimes_R B_y) \xrightarrow{\sim}H^{i}(B_x\otimes_R B_y)
\end{equation} for all $i$.
\end{theorem}

For an arbitrary reflection $s\in W$, this theorem yields an isomorphism $B_s(-1)\simeq H^{-1}(B_s\otimes_R B_s)\simeq H^1(B_s\otimes_R B_s)\simeq B_s(+1)$. Using canonical projection to the lowest graded summand and canonical inclusion from the highest graded summand, we can define the \emph{asymptotic Hecke category}. More generally, an object lying over $j_z$ (for a summand of $j_xj_y$) should correspond to both the lowest and highest graded parts of $B_xB_y$. Relative hard Lefschetz identifies these via $\eta^i$. The maps corresponding to the tensor product take the form:

\begin{equation}\label{eq:inc_proj_mon}
	\begin{tikzcd}[column sep=40pt]
		H^{-i}(B_xB_y) \ar[r,dashed,"{\eta^i}"] \ar[d,hookrightarrow,swap,"{\text{canonical inclusion}}"]& H^{i}(B_xB_y) \\
		B_xB_y\ar[r,swap,"{\eta^i}"] & B_xB_y\ar[u,twoheadrightarrow,swap,"{\text{canonical projection}}"] \\
	\end{tikzcd}
\end{equation}

This motivates the definition of a monoidal category categorifying $J_c$.  The following is a compression of the construction of \cite[Section 5]{elias.will16}.

\begin{construction}
	Fix a two-sided Kazhdan--Lusztig cell $c$ with $a$-value $i$ in a Coxeter system $(W,S)$. 
	\begin{itemize}
		\item We define the subcategory $\mathcal{H}_W^{<c}$ as the full subcategory of $\mathcal{H}_W$ generated by objects $B_x$ such that $x<_Jc$. Let $\mathcal{I}^c$ denote the tensor ideal of morphisms in $\mathcal{H}_W$ factoring over objects of $\mathcal{H}_W^{<c}$. We define the quotient category by $(\mathcal{H}_W^{c})'\coloneqq \mathcal{H}_W/\mathcal{I}^c$.
		
		\item Inside $(\mathcal{H}_W^c)'$, we restrict to the full graded additive subcategory $\tilde{\mathcal{H}}_W^c$ generated only by objects $B_x$ for $x\in c$. 
		
		\item We now enrich the grading-free full subcategory $\mathcal{H}_W^c$ of $\tilde{\mathcal{H}}_W^c$ (i.e. the subcategory generated only by $B_x$ without shifts) with a new monoidal product using the $i$-th perverse cohomology: \begin{equation}\label{eq:asymptotic_monoidal_product}
			B_x\star B_y \coloneqq H^{-i}(B_xB_y) \in \mathcal{H}_W^c.
		\end{equation} 

	\end{itemize}
\end{construction}

\begin{remark}\label{rem:explanation_quotient}
	The quotient construction of $(\mathcal{H}_W^c)'$ is necessary to account for the fact that in the construction of the $J$-ring one discards any summand of $j_xj_y$ lying in lower cells. The perverse filtration of $\mathcal{H}_W$ descends to $\tilde{\mathcal{H}}_W^c$ and $\mathcal{H}_W^c$ and any $H^{-i}(B_xB_y)$ contains no summands of lower cells.
	
	For the monoidal structure on $\mathcal{H}_W^c$ we use inclusions and projections as in \eqref{eq:inc_proj_mon}. Furthermore, we note that the category is rigid and pivotal using the canonical isomorphism of a Soergel bimodule to its double dual by \cite[Proposition 5.6(3)]{elias.will16}. This is always implicitly used, for example for the trace of \Autoref{def:s-matrix}.
	
	In all cases we consider the pivotal structure is also spherical. For the dihedral group we get this for $\mathcal{H}^c_W$ as the definition of $\mathcal{H}_W$ is symmetric in both generators. In all other cases, we look at small diagonal cells such that each object in $\mathcal{H}^h_W$ is self-dual, hence it is spherical because we are pivotal.
\end{remark}

\section{The center of the asymptotic Hecke category}\label{sec:center_asymptotic}
We follow the categorical notation of \cite{EGNO}. After recalling the definition of multifusion categories and properties of their Drinfeld centers, we show that studying the asymptotic Hecke category of a finite $J$-cell $c$ reduces to studying a considerably smaller $H$-cell $h\subset c$.

\subsection{Multifusion categories}

Let $\Bbbk$ be an algebraically closed field. Outside this section we assume that all categories we consider are $\Bbbk=\mathbb{C}$-linear.

\begin{definition}[{\cite[Section 4.1]{EGNO}}]
	A category $\mathcal{C}$ is \emph{multifusion} if it is locally finite, $\Bbbk$-linear, abelian, rigid, monoidal, and semisimple, with bilinear tensor product $\otimes: \mathcal{C}\times \mathcal{C}\to \mathcal{C}$ and finitely many simple objects. If, furthermore, $\End_{\mathcal{C}}(\mathds{1})\simeq \Bbbk$ for the monoidal unit $\mathds{1}$, we call $\mathcal{C}$ a \emph{fusion category}.
\end{definition}

\begin{example}
	Examples of fusion categories include $\mathrm{Vec}(G)\coloneqq \mathrm{Vec}_{\Bbbk}(G)$, the category of $G$-graded finite-dimensional $\Bbbk$-vector spaces, and $\mathrm{Rep}(G)\coloneqq \mathrm{Rep}_{\Bbbk}(G)$, the category of representations of a finite group $G$ over $\Bbbk$ when $\mathrm{char}(\Bbbk)$ and $|G|$ are coprime.
\end{example}

\begin{remark}\label{rem:multifusion_asymptotic}
	Let $K(\mathcal{C})$ denote the Grothendieck ring of a multifusion category $\mathcal{C}$. By definition, it is a multifusion ring by choosing the equivalence classes of the simple objects as basis elements. 
	
	By \cite[Section 5.2]{elias.will16}, the asymptotic Hecke category is rigid and pivotal. We have seen in \Autoref{rem:duflo} that the asymptotic Hecke algebra $J_c$ is multifusion if $c$ is finite. Therefore, $\mathcal{H}_W^c$ is a multifusion category. The sum $\bigoplus_{d\in D_c}B_d$ for $D_c$ the intersection of the Duflo involutions $D$ with the cell $c$ is then the unit of $\mathcal{H}_W^c$.
\end{remark}

By \cite[Theorem 4.3.1]{EGNO}, in a multifusion category $\mathcal{C}$ the space $\mathrm{End}_{\mathcal{C}}(\mathds{1})$ is semisimple, hence $\mathds{1}=\bigoplus_{i\in I}\mathds{1}_i$ for non-isomorphic indecomposable objects~$\mathds{1}_i$.

\begin{definition}
	Let $\mathcal{C}$ be a multifusion category with $\mathds{1}=\bigoplus_{i\in I}\mathds{1}_i$. For $i,j\in I$, the \emph{component subcategory} $\mathcal{C}_{ij}\coloneqq \mathds{1}_i\otimes \mathcal{C}\otimes \mathds{1}_j$ is the full subcategory generated by objects $\mathds{1}_i\otimes X \otimes \mathds{1}_j$.
\end{definition}

As abelian categories, $\mathcal{C} \simeq \bigoplus_{i,j\in I}\mathcal{C}_{i,j}$. The monoidal product maps $\mathcal{C}_{ij}\times \mathcal{C}_{jk}\to \mathcal{C}_{ik}$, and duals of objects in $\mathcal{C}_{ij}$ lie in $\mathcal{C}_{ji}$ \cite[Remark 4.3.4]{EGNO}. Each $\mathcal{C}_{ii}$ is a fusion category with unit $\mathds{1}_i$. The Drinfeld center of a multifusion category is equivalent to that of its fusion subcategories, as we show next and apply to the asymptotic Hecke category in \Autoref{sub:h-reduction}.

\subsection{The Drinfeld center of multifusion categories} 

We recall the definition of the Drinfeld center, see \cite[Definition 7.13.1]{EGNO}.

\begin{definition}
	Let $\mathcal{C}$ be a monoidal category. The \emph{center} $\mathcal{Z}(\mathcal{C})$ is a category with objects $(Z,\gamma)$ where $Z\in\mathcal{C}$ and $\gamma$ is a family of natural morphisms $\gamma_X:X\otimes Z \to Z \otimes X$ for all $X\in \mathcal{C}$ satisfying the hexagon axiom.
\end{definition}

Most properties of $\mathcal{C}$ transfer to $\mathcal{Z}(\mathcal{C})$. For example, the center is always a monoidal category, and it is also fusion if $\mathcal{C}$ is, see \cite[Theorem 9.3.2]{EGNO}. The Drinfeld center $\mathcal{Z}(\mathcal{C})$ is a special case of a \emph{braided monoidal category}, i.e., a monoidal category $\mathcal{D}$ equipped with a natural family of isomorphisms $c_{X,Y}: X \otimes Y \to Y \otimes X$
for all objects $X,Y \in \mathcal{D}$ satisfying the hexagon axioms. The braiding on $\mathcal{Z}(\mathcal{C})$ arises from the \emph{half-braidings} $\gamma_X$ defining its objects.

If the center admits a spherical structure, i.e., a pivotal structure such that the dimension of every object agrees with that of its dual, then the $S$-matrix can be defined. The following definition extends naturally to spherical braided categories.

\begin{definition}\label{def:s-matrix}
	Let $\mathcal{C}$ be a spherical fusion category and $\mathcal{Z}(\mathcal{C})$ its Drinfeld center. For $(Z_i,\gamma^i)_i$ a complete list of all simple objects of $\mathcal{Z}(\mathcal{C})$ we define the \emph{$S$-matrix} of $\mathcal{Z}(\mathcal{C})$ to be \begin{equation}
		S\coloneqq (\mathrm{tr}(\gamma^i_{Z_j}\circ \gamma^j_{Z_i}))_{i,j},
	\end{equation} where $\mathrm{tr}$ denotes the trace of an endomorphism $f:Z\to Z$, i.e. the element in $\Bbbk$ corresponding to $f$ after applying the evaluation and coevaluation, see \cite[Section 8.13]{EGNO}.
\end{definition}

\begin{example}\label{ex:center_vecg}
	Let $G$ be a finite group. The Drinfeld center of $\mathrm{Vec}(G)$ is completely described, see \cite[Example 4.15.4]{EGNO}. It is $\mathcal{Z}(\mathrm{Vec}(G))\simeq (\mathrm{Vec}(G))^G$, the category of $G$-equivariant $G$-graded vector spaces where $G$ acts on $\mathrm{Vec}(G)$ by conjugation. The simple objects are in correspondence with the set of pairs $(C,V)$, where $C$ is a conjugacy class of $G$ and $V$ is a simple representation of the centralizer of some $g\in C$. Duals are given by inverting on the group level and we are trivially spherical.
	
	For $G=S_3$ this gives for example eight simple objects in the center, three lying over the trivial conjugacy class, three over the conjugacy class of the $3$-cycle and two over the conjugacy class of the $2$-cycle. By \cite[Example 8.13.6]{EGNO} the $S$-matrix is \begin{equation}\label{eq:smatrix}
		S_{(C,V),(C',V')} = \frac{\abs{G}}{\abs{C_G(a)}\abs{C_G(a')}}\sum_{g\in G(a,a')}\mathrm{tr}_V(ga'g^{-1})\mathrm{tr}_{V'}(g^{-1}ag),
	\end{equation} where $a\in C,a'\in C'$ and $G(a,a')=\{g\in G\mid aga'g^{-1}=ga'g^{-1}a\}$.
\end{example}

For a multifusion category $\mathcal{C}$, the center reduces to the centers of fusion subcategories $\mathcal{C}_{ii}$ for $i\in I$. We call $\mathcal{C}$ \emph{indecomposable} if $I$ cannot be partitioned into non-empty subsets $I=J\coprod K$ with $\mathcal{C}_{j,k}=0$ for all $j\in J$ and $k\in K$. For decomposable $\mathcal{C}$, the center is the direct sum of centers of summands. For indecomposable $\mathcal{C}$, the center equals that of any fusion subcategory $\mathcal{C}_{ii}$.

\begin{proposition}[{\cite[Theorem 2.5.1]{KONG2018749}}]\label{prop:centr_reduction}
	For an indecomposable multifusion category $\mathcal{C}$ with component fusion subcategories $\mathcal{C}_{ii}$ for $1\leq i \leq n$ we have \begin{equation}
		\mathcal{Z}(\mathcal{C}) \simeq \mathcal{Z}(\mathcal{C}_{ii}).
	\end{equation}
	Thus, the center of an indecomposable multifusion category is fusion.
\end{proposition}

We apply this result in the next section to the asymptotic Hecke category.

 \subsection{$H$-cell reduction}\label{sub:h-reduction}
 
 We reduce the computation of the center of the asymptotic Hecke category associated to a finite $J$-cell to that of an $H$-cell contained inside. This process, called \emph{$H$-reduction} or Clifford--Munn--Ponizovski\u\i{} theory has been applied to monoids, algebras and categories, see for example \cite[Theorem 15]{mackaay20}.
 
 Let $W$ be a Coxeter group and $c\subset W$ a finite $J$-cell. The decomposition of the asymptotic Hecke category $\mathcal{H}_W^c$ into component subcategories comes from the decomposition of $c$ into left and right cells. We write $c^L_i$ and $c^R_i$, where $1\leq i \leq n$, for a list of left and right cells. By \Autoref{rem:duflo,rem:directsum} the monoidal unit is the direct sum of all objects lying over Duflo involutions in the $J$-ring: $\mathds{1}_{\mathcal{H}_W^c}=\bigoplus_{1\leq i \leq n}B_{d_i}$, where $\{d_i\}$ is the set of all Duflo involutions of $W$ lying in $c$. We assume that the left and right cells are ordered such that $d_i\in c^L_i\cap c^R_i$.
 
 \begin{definition}
 	We call a non-empty intersection of a left and a right cell  an \emph{$H$-cell}. If an $H$-cell contains a Duflo involution we call it \emph{diagonal}.
 \end{definition}
 
 Any diagonal $H$-cell with $d_i\in h_i=c^L_i\cap c^R_i\subset W$ gives a component subcategory $\mathcal{H}_W^h\coloneqq (\mathcal{H}_W^c)_{ii}=B_{d_i}\star \mathcal{H}_W^c\star B_{d_i}$ of $\mathcal{H}_W^c$. This is a fusion category and we have \begin{equation} \label{equ:h_reduction_center}
 	\mathcal{Z}(\mathcal{H}_W^h) \simeq \mathcal{Z}(\mathcal{H}_W^c)
 \end{equation} by \Autoref{prop:centr_reduction}.
 
 Hence, the computation of the Drinfeld center of the asymptotic Hecke category of a $J$-cell reduces to that of an $H$-cell. 

\subsection{The centers of the asymptotic Hecke category for finite Weyl groups}\label{sub:center_weyl}

For finite Weyl groups, we identify the asymptotic Hecke categories -- viewed as fusion categories -- with known categories arising from classical geometric results. We give an overview on the classification and describe their centers and $S$-matrices using $H$-reduction.

By \cite[Chapter 4]{lusz1984charac} we have an assignment of a two-sided Kazhdan--Lusztig cell $c$ in a Weyl group to a finite group $G_c$ and an embedding $c\to M(G_c)$, where $M(G_c)$ consists of tuples $(g,V)$ for $g\in G_c$ unique up to conjugacy and $V$ a simple representation of the centralizer of $g$. 

For any left cell $c^L\subseteq c$ there is further an association to a subgroup $H_{c^L}\leq G_c$ in \cite{lusztig1987Leading} such that the asymptotic Hecke algebra $J_h$ associated to the $H$-cell $h\coloneqq c^L \cap (c^L)^{-1}$ is, as a based (or multifusion) ring, conjectured to be isomorphic to $K_{G_c}(G_c/H_{c^L} \times G_c/H_{c^L})$, which is short for the Grothendieck ring of $\mathrm{Coh}_{G_c}(G_c/H_{c^L}\times G_c/H_{c^L})$, the category of $G_c$-equivariant coherent sheaves on the set $(G_c/H_{c^L})^2$. Furthermore, the conjecture \cite[Conjecture 3.15]{lusztig1987Leading}, is extended to the claim that the disjoint union $X\coloneqq \coprod_{c^L\subset c}{G_c/H_{c^L}}$ gives a multifusion ring isomorphic to $K_{G_c}(X\times X)\simeq J_c$. This was proven by Lusztig himself in the case that $G_c$ is abelian. A complete proof was achieved by Bezrukavnikov, Finkelberg and Ostrik in \cite[Theorem 4]{Bez2009}. For all but three exceptions in type $E_7$ and $E_8$, they even showed that $J_c$ is categorified by $\mathrm{Coh}_{G_c}(X\times X)$ for the same $G_c$-set $X$. We call the three exceptions the \emph{exceptional} cells.

The results presented above are summarized in the following example:

\begin{example} \label{ex:classification}
	The categories $\mathcal{H}_W^h$ for a diagonal $H$-cell $h=c^L\cap (c^L)^{-1}$ in a non-exceptional two-sided Kazhdan--Lusztig cell $c$ of a finite Weyl group $W$ are given by $\mathrm{Coh}_{G_c}(G_c/H_{c^L}\times G_c/H_{c^L})$, i.e. categories of equivariant coherent sheaves on a finite set, for the following possibilities of $G_c$ and $H_{c^L}$. Here $S_n$ denotes the symmetric group on $n$ letters, $D_m$ stands for the dihedral group of order $2m$.
	\begin{itemize}
		\item In type $A_n$ any $H$-cell has size one, we always have $G_c=\{\star\}=H_{c^L}$.
		\item In type $B_n$ the size of an $H$-cell is $2^k$ for some $k^2+k\leq n$. The groups $G_c$ and $H_{c^L}$ are some elementary abelian $2$-groups.
		\item In type $D_n$ we have the same result as in $B_n$ except that $k^2\leq n$.
		\item In type $E_6$ to $E_8$ the group $G_c$ is a symmetric group on at most five letters $S_1,\ldots, S_5$: \begin{itemize}
			\item For $G_c=S_3$ we can have $H_{c^L}\in \{S_1,S_2,S_3\}$.
			\item  For $G_c=S_4$ we can have $H_{c^L}\in \{S_2,S_2\times S_2, S_3, D_4, S_4\}$.
			\item  For $G_c=S_5$ we can have $H_{c^L}\in \{S_2,S_2\times S_2,S_3,D_4,S_2\times S_3,S_4,S_5\}$.
		\end{itemize}
		\item In type $F_4$ we get $G_c < S_4$ with the same possible subgroups for $H_{c^L}$ as before.
		\item In type $G_2$ we get $G_c\in \{S_1,S_3\}$, where for $G_c=S_3$ only $H_{c^L}=S_2$ occurs. 
	\end{itemize}
\end{example}

\begin{remark}\label{rem_center_hreduction}
	We want to motivate the connection of the set $M(G_c)$ to the center of $\mathrm{Coh}_{G_c}(X\times X)$. The categories $\mathcal{C}\coloneqq \mathrm{Coh}_{G_c}(X\times X)$ are multifusion. If $X=\cup X_i$ is a disjoint union of transitive $G_c$-sets $X_i$, the categories $\mathcal{C}_{ij}\coloneqq \mathrm{Coh}_{G_c}(X_i\times X_j)$ are component subcategories. By \Autoref{prop:centr_reduction} the centers of $\mathcal{C}_{ii}$ and $\mathcal{C}$ are equivalent. 
	
	If one chooses $X=G_c$ we have $\mathrm{Coh}_{G_c}(X\times X)\simeq \mathrm{Vec}(G_c)$. Therefore, the center $\mathcal{Z}(\mathcal{C})$ is equivalent to the center of the category of $G_c$-graded vector spaces. Indeed, the set $M(G_c)$ has an analogous description as the simple objects of the center $\mathcal{Z}(\mathrm{Vec}(G_c))\simeq (\mathrm{Vec}_{G_c})^{G_c}$ as seen in \Autoref{ex:center_vecg}. On the one hand we have tuples $(g,V)$ of an element $g\in G_c$ up to conjugacy and a representation of its centralizer; and on the other hand we have tuples $(C,V)$ of the conjugacy class and representations of its stabilizer.

	Furthermore, the $S$-matrix computed for $\mathcal{Z}(\mathrm{Vec}(G))$ is very similar to the pairing on $M(G_c)$ defined in \cite[Equation 4.14.3]{lusz1984charac}: \begin{equation}\label{eq:s_matrix}
		\{(x,\sigma),(y,\tau)\} \coloneqq \sum_{g\in G_c, xgyg^{-1}=gyg^{-1}x}\frac{\mathrm{tr}(g^{-1}x^{-1}g,\tau)\mathrm{tr}(gyg^{-1},\sigma)}{\abs{C_{G_c}(x)}\abs{C_{G_c}(y)}}.
	\end{equation}
	These are nearly the entries of the \emph{normalized} $S$-matrix of the Drinfeld center. By normalized, we mean that we divided by $\abs{G}$, the square root of the categorical dimension of  $\mathcal{Z}(\mathrm{Vec}(G_c))$; see \cite[Section 8.14]{EGNO}. 
	
	The only difference is the inverse of $x$, which permutes the entries of the matrix. While the square of Lusztig's matrix is the identity, the square of \eqref{eq:smatrix} returns $\abs{G}^2$ times a permutation matrix describing duality of the simple objects. This explanation can be found in the \href{https://webusers.imj-prg.fr/~jean.michel/gap3/htm/chap098.htm\#I21}{GAP3-Documentation, see Section 98.22 `DrinfeldDouble'}, from Jean Michel.\footnote{\url{https://webusers.imj-prg.fr/~jean.michel/gap3/htm/chap098.htm\#I21}, last updated November 27, 2003, accessed on October 6, 2025.}
	
	In all cases $G_c$ we consider, the representation are self-dual, hence for us the normalized $S$-matrix and Lusztig's matrix of pairings are the same.
\end{remark}

With the description of $\mathcal{H}^h_W$ from \Autoref{ex:classification} and the $S$-matrix data from \eqref{eq:s_matrix} of \Autoref{rem_center_hreduction} together with \Autoref{prop:centr_reduction} we get the following $S$-matrices for the two-sided cells in Weyl groups:

\begin{corollary}\label{cor_smatrix_weyl}
	Let $c$ be a non-exceptional two-sided Kazhdan-Lusztig cell in a finite Weyl group $W$. The asymptotic Hecke category associated to $c$ as well as the $S$-matrix of its center is one of the following cases:\begin{itemize}
		\item For any $c$ where a diagonal $H$-cell has size one we have $\mathcal{H}_W^c=\mathrm{Coh}(X\times X)$ where $X$ has the same cardinality as the number of left and right cells in $c$. We have $\mathcal{H}_W^h\simeq \mathrm{Coh}(\star)\simeq \mathrm{Vec}$ for any diagonal $H$-cell. The center $\mathcal{Z}(\mathcal{H}_W^c)\simeq \mathcal{Z}(\mathcal{H}_W^h)\simeq \mathrm{Vec}$ has size one and the $S$-matrix is \begin{equation}
			S_c=\begin{pmatrix}
				1\\
			\end{pmatrix}.
		\end{equation}  This happens for any cell in type $A_n$ and also for all cells containing only the trivial element. More examples of cells can be found in \cite[Section 8]{Mackaay2023}.
		
		\item If the asymptotic Hecke category of $c$ is equivalent to $\mathrm{Coh}_G(X\times X)$ for an elementary abelian $2$-group, i.e. $G_c\simeq (\ZZ/2\ZZ)^k$, we have $\mathcal{Z}(\mathcal{H}_W^c)\simeq \mathcal{Z}(\mathrm{Vec}(G_c))\simeq \bigboxtimes_{1\leq i \leq k} \mathcal{Z}(\mathrm{Vec}(\ZZ/2\ZZ))$, i.e. the center is a $k$-fold \emph{Deligne tensor product}. The center then contains $4^k$ simple objects and the $S$-matrix is the $k$-fold Kronecker product of the $S$-matrix of $\mathcal{Z}(\mathrm{Vec}(\ZZ/2\ZZ))$, which is \begin{equation}
			S(\mathcal{Z}(\mathrm{Vec}(\ZZ/2\ZZ)))=\begin{pmatrix}
				1 & 1 & 1 & 1 \\
				1 & 1 & -1 & -1 \\
				1 & -1 & 1 & -1 \\
				1 & -1 & -1 & 1 \\
			\end{pmatrix}.
		\end{equation} Since the dimension of $\mathrm{Vec}(\ZZ/2\ZZ)$ is $2$, the normalization agrees with the table in \cite[Section 4.15]{lusz1984charac}.
	
		\item If the asymptotic Hecke category of $c$ is equivalent to $\mathrm{Coh}_G(X\times X)$ with $G=S_3$ the center of the asymptotic Hecke category is $\mathcal{Z}(\mathrm{Vec}(S_3))$ which has eight simple objects and the $S$-matrix is \begin{equation}
			S(\mathcal{Z}(\mathrm{Vec}(S_3))) = \begin{pmatrix}
				4 & 2 & 2 & 0 & 0 & -2 & -2 & 2 \\
				2 & 1 & 1 & -3 & -3 & 2 & 2 & 2 \\
				2 & 1 & 1 & 3 & 3 & 2 & 2 & 2 \\
				0 & -3 & 3 & 3 & -3 & 0 & 0 & 0 \\
				0 & -3 & 3 & -3 & 3 & 0 & 0 & 0 \\
				-2 & 2 & 2 & 0 & 0 & 4 & -2 & -2 \\
				-2 & 2 & 2 & 0 & 0 & -2 & -2 & 4 \\
				-2 & 2 & 2 & 0 & 0 & -2 & 4 & -2 \\
			\end{pmatrix}.
		\end{equation} Normalization by the dimension of $ \dim(\mathrm{Vec}(S_3))=6$ gives the table of \cite[Section 4.15]{lusz1984charac}. Note that some rows have been left out in that source, they are permutations of some rows given.
		\item If the asymptotic Hecke category of $c$ is equivalent to $\mathrm{Coh}_G(X\times X)$ with $G=S_4$ the center of the asymptotic Hecke category is $\mathcal{Z}(\mathrm{Vec}(S_4))$ which has $21$ simple objects. To count this we have to compute all centralizer subgroups and count their irreducible representations. The matrix can also be found in \cite[Section 4.15]{lusz1984charac}.
		
		\item If the asymptotic Hecke category of $c$ is equivalent to $\mathrm{Coh}_G(X\times X)$ with $G=S_5$ the center of the asymptotic Hecke category is $\mathcal{Z}(\mathrm{Vec}(S_5))$ which has $39$ simple objects, again see \cite[Section 4.15]{lusz1984charac}.
	\end{itemize}
\end{corollary}

\subsection{The exceptional cells in Weyl groups}\label{sub:exceptional}

In the three exceptional cases in type $E_7$ and $E_8$ we have a categorification of $\mathcal{H}_W^c$ by \cite[Theorem 1.1]{ost2014tensor}: 

\begin{theorem}
	For an exceptional cell $c$ in type $E_7$ or $E_8$, there is a tensor equivalence $\mathcal{H}_W^c\simeq \mathrm{Vec}^\omega(\ZZ/2\ZZ)\boxtimes \mathrm{Coh}(Y'\times Y')$. 
\end{theorem}

Note, that the category $\mathcal{H}_W^c$ is denoted by $\mathcal{P}_c$ in \cite{ost2014tensor}. The set $Y'$ has cardinality $512$ for the exceptional cell in type $E_7$ and $4096$ for the two exceptional cells in type $E_8$. The cardinality of the set $Y'$ gives the number of left or right cells in $c$, the $H$-cells have only size two and are therefore categorified by $\mathrm{Vec}^\omega(\ZZ/2\ZZ)$, where $\omega$ denotes the non-trivial twist.

\begin{corollary}\label{cor_e78}	
	Let $c\subset W$ be an exceptional cell in type $E_7$ or $E_8$. The center of the asymptotic Hecke category associated to $c$ is $\mathcal{Z}(\mathcal{H}^c)\simeq \mathcal{Z}(\mathrm{Vec}^\omega(\ZZ/2\ZZ))$ for $\omega$ a non-trivial twist. We have four simple objects in $\mathcal{Z}(\mathcal{H}^c_W)$ and the $S$-matrix is \begin{equation}
		S(\mathcal{Z}(\mathrm{Vec}^\omega(\ZZ/2\ZZ)))=\begin{pmatrix}
			1 & 1 & 1 & 1 \\
			1 & 1 & -1 & -1 \\
			1 & -1 & -1 & 1 \\
			1 & -1 & 1 & -1 \\
		\end{pmatrix}. 
	\end{equation}
\end{corollary}

\section{The dihedral case}\label{sec:dihedral}

We give a complete description of the asymptotic Hecke category associated to a dihedral group. We will see that the category is known in the literature as the even or adjoint part of the Verlinde category. The Drinfeld center of the Verlinde category and its adjoint are also known. We will combine all the different known results to then prove Lusztig's conjecture for dihedral groups. 

Let $W=\langle s,t \mid s^2=t^2=(st)^p=1 \rangle$ be the Coxeter group of type $I_2(p)$ where $p \geq 3$.  

\subsection{The asymptotic Hecke algebra for dihedral groups}

All data on $h_{x,y,z}$ and the asymptotic Hecke algebra are known, see for example \cite[Section 4]{ducloux2006731}. There are always three two-sided cells for $p\geq 3$. 

The unit always forms its own two-sided cell $c_0=\{1\}$ as $x\leq_K 1$ for all $x\in W$ and $K\in \{L,R,J\}$ since $b_x=b_1b_xb_1$. The $a$-value is $0$. Similarly, the longest word $c_p=\{w_0\}$ for \begin{equation}
w_0=\underbrace{sts\ldots}_{p \text{ times}}
\end{equation} forms its own two-sided cell with $a$-value $p$ as $b_xb_{w_0}=q(v)b_{w_0}$ for some polynomial $q(v)\in \mathbb{Z}[v^{\pm 1}]$. Furthermore, we have the so-called \emph{subregular} cell of $a$-value $1$. It contains any non-trivial word that has a unique reduced expression. These are all remaining elements $c_1=\{s,st,sts,\ldots,t,ts,tst,\ldots\}$. The left and right cells are characterized by the right and left descending sets. We can visualize the cell structure in a box diagram, where the big boxes correspond to $J$-cells, columns to $R$-cells, rows to $L$-cells and small boxes to $H$-cells.

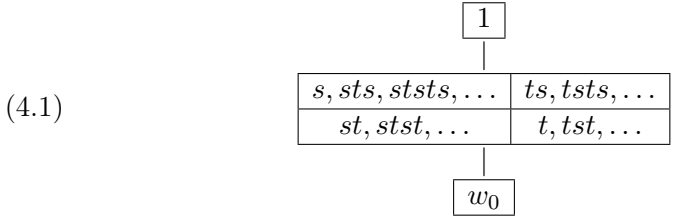
\begin{equation}
	\begin{tabular}{c}			
		\begin{tabular}{|c|}
			\hline
			$	1$	 \\
			\hline
		\end{tabular}\\			
		| \\
		\begin{tabular}{|c|c|}
			\hline
			$	s, sts, ststs,\ldots$ & $ts, tsts,\ldots$ \\
			\hline
			$	st, stst, \ldots$	& $t,tst,\ldots$ \\
			\hline
		\end{tabular}
		\\
		| \\	
		\begin{tabular}{|c|}
			\hline
			$w_0$\\
			\hline
		\end{tabular}		
	\end{tabular}
\end{equation} 

The multiplication table of the $J$-ring can also be found in \cite[Section 4]{ducloux2006731}. The coefficients $\gamma_{x,y,z}$ are either $0$ or $1$. We denote by $s_k$ the unique word of length $k$ starting in $s$ and by $t_l$ the unique word of length $l$ starting with $t$ for $k,l<p$. The multiplication in the $J$-ring is then:
\begin{equation}\label{eq:mult_jring}
	j_{s_k} j_{a_l} =
	\begin{cases}
		0 & \text{if } 
		\begin{aligned}
			&k \text{ even and } a = s \\[-2pt]
			&\text{or } k \text{ odd and } a = t
		\end{aligned} \\[3pt]
		\sum_{u=\max\{0,k+l-p\}}^{\min\{k,l\}-1} j_{s_{k+l-1-2u}}
		& \text{otherwise.}
	\end{cases}
\end{equation}

In type $I_2(5)$ this is for example:
\begin{equation}\label{table:muli24}
	\begin{tabular}{c|cccc}
		$\cdot$ & $j_s$ & $j_{st}$ & $j_{sts}$ & $j_{stst}$ \\
		\hline
		$j_s$ & $j_s$ & $j_{st}$ & $j_{sts}$ & $j_{stst}$ \\
		$j_{ts}$ & $j_{ts}$ & $j_{t} + j_{tst}$ & $j_{ts}+ j_{tsts}$ & $j_{tst}$\\
		$j_{sts}$ & $j_{sts}$ & $j_{st} + j_{stst}$ & $j_s + j_{sts}$ & $j_{st}$\\
		$j_{tsts}$ & $j_{tsts}$ & $j_{tst}$ & $j_{ts}$ & $j_t$\\
	\end{tabular}
\end{equation}
We can directly read off that the unit is $j_s+j_t$. In the next subsection, we observe that categories with such Clebsch--Gordan-like multiplication rules have been extensively studied, which allows us to find all possible categorifications of the $J$-ring.

\subsection{$A_n$-fusion rule categories}\label{sub:type_An}

\begin{definition}
	We say that a fusion category $\mathcal{C}_n$ has \emph{$A_n$-fusion rules } if it has $n$ simple objects, which we may label by $X_0,\ldots,X_{n-1}$, such that the fusion graph showing the monoidal product by $X_1$ is the Dynkin diagram of type $A_n$. 
	
	This means $X_1\otimes X_0\simeq X_1\simeq X_0\otimes X_1$, $X_{1}\otimes X_{n-1}\simeq X_{n-2}\simeq X_{n-1}\otimes X_1$ and $X_1\otimes X_i\simeq X_{i-1}\oplus X_{i+1}$ for all $1\leq i \leq n-2$.
\end{definition}

\begin{remark}
	The fusion data of these categories are mentioned in \cite[Section 2.2]{lus94-exotic} under the name \emph{Verlinde--Wess--Zumino--Witten}. A complete classification of all possible categorification is done in \cite[Prop 8.2.3]{FroehlichKerler1993}. There, they call it the semisimple quotient of the representation category of $U_q(\mathfrak{sl}_2)$.
\end{remark}

An overview of the categorical data of fusion categories with $A_n$-fusion rules  can be found in \cite{edie-morrison-17}. All associators have been classified by \cite{FroehlichKerler1993}. It turns out that the evaluation of the loop is an invariant of the category and there are only finitely many such values that will define a category. 

\begin{lemma}{\cite[Proposition 8.2.3]{FroehlichKerler1993}}\label{lem:classification_ak}
	For every $1\leq l \leq n$ coprime to $n+1$ there exists a fusion category, unique up to natural equivalence, with $A_{n}$-fusion rules, such that the composition of coevaluation and evaluation (compare to \cite[Section 2.10]{EGNO}) is: $\mathrm{coev}_{X_1}\circ \mathrm{ev}_{X_1}=2\cos\left(\frac{l\pi}{n+1}\right)$. We denote this category by $\mathcal{C}_n^l$.
\end{lemma}

To work with these categories one introduces the notion of quantum numbers. Let $q\coloneqq e^{\frac{\pi i}{n+1}}$ be a $2(n+1)$-th root of unity, then we set the quantum numbers, depending on $l$, to be $[0]\coloneqq 0, [1]\coloneqq 1$, $[2]\coloneqq q^l+q^{-l}$, and inductively $[k]\coloneqq [2][k-1]-[k-2]$. Hence, the invariant of \Autoref{lem:classification_ak} is $[2]$. 

\begin{example}\label{ex:twistedZZ2ZZ}
	There are detailed formulas for all associators found in \cite{edie-morrison-17}, the calculations have been done in \autocite[Chapter 9]{kauff1994}. 
	
	For $n=2$, the monoidal product on a category with $A_2$-fusion rules is 
	\begin{equation}
		X_0\otimes X_0 \simeq X_0 \simeq X_1 \otimes X_1,\quad X_0\otimes X_1 \simeq X_1 \simeq X_1 \otimes X_0.
	\end{equation}
	
	The associators are mostly trivial. Only on the triple $(X_1,X_1,X_1)$ we get that $\alpha_{X_1,X_1,X_1}: X_1 \simeq (X_1\otimes X_1) \otimes X_1 \to X_1 \otimes (X_1\otimes X_1) \simeq X_1$ is $-\frac{[1]}{[2]}$ times the identity.  For $l=1$ we have $[2] = 1$, giving us the category $\mathcal{C}_2^1 \simeq \mathrm{Vec}^\omega(\mathbb{Z}/2\mathbb{Z})$ of twisted $\mathbb{Z}/2\mathbb{Z}$-graded vector spaces, while for $l=2$ we get $[2]=-1$ and the non-twisted version.
\end{example}

\subsection{The adjoint part of $A_n$-fusion rule categories}\label{sub:adjoint_category}

In \cite{ediemichell2018brauerpicard} the subcategory of $\mathcal{C}_n$ generated by the even elements is called the \emph{adjoint subcategory}. An explanation for this term can be found in \cite[Section 3.6 and 4.14]{EGNO}. 

\begin{remark}
	For a based ring $A$ with basis $B=\{b_i\}$ we call the smallest subring $A_{ad}\subset A$, such that all $b_ib_i^*$ lie in $A_{ad}$ the \emph{adjoint subring}. For a fusion category $\mathcal{C}$ we write $\mathrm{Ad}(\mathcal{C})$ for the full fusion subcategory such that $K(\mathrm{Ad}(\mathcal{C}))=K(\mathcal{C})_{ad}$ and call it the \emph{adjoint subcategory}.
	
	In $\mathcal{C}_n$ all objects are self-dual, and any monoidal product $X_i\otimes X_i$ decomposes into a sum of even summands $X_{2j}$. This comes from the fact that $\mathcal{C}_n$ is universally $\ZZ/2\ZZ$-graded in the sense of \cite[Section 4.14]{EGNO} and the adjoint part is the trivial component of the universal grading on $\mathcal{C}_n$.
\end{remark}

While we saw in \Autoref{lem:classification_ak} that the categories $\mathcal{C}_n^l$ are the only categorifications of Verlinde-type fusion rings, it is not yet clear that the adjoint subcategories are the only possibilities for categorifications of the adjoint fusion rings. Recent work by Etingof and Ostrik \cite{etingof2022semisimplification} shows that this is indeed the case. In their shorthand notation we write $K_n$ for the Grothendieck ring of the adjoint part of $\mathcal{C}_{2n+1}$ and $K'_n$ for the Grothendieck ring of the adjoint part of $\mathcal{C}_{2n+2}$. 

\begin{lemma}\label{rem:ostriksresult}
	Let $\mathcal{C}$ be a pivotal fusion category categorifying the fusion ring $K_n$ or $K'_{n'}$ for $n> 2$ or $n'\geq 0$. Then there is a tensor equivalence $\mathcal{C}\simeq \mathrm{Rep}(\mathfrak{so}(3)_q)$ for $q$ a primitive $4(n+1)$-th root of unity, respectively a primitive $4(n'+2)$-th root of unity.
\end{lemma}

		This is \cite[Theorem A.3 and Remark A.4(ii)]{etingof2022semisimplification}. Here $\mathrm{Rep}(\mathfrak{so}(3)_q)$ is the fusion category of tilting modules over the quantum enveloping algebra of $\mathfrak{so}(3)$ specialized at the root of unity $q$. In our notation these are the categories $\mathrm{Ad}(\mathcal{C}_{2n+1}^l)$, respectively $\mathrm{Ad}(\mathcal{C}_{2n+2}^l)$, for the right choice of $l$.

\begin{remark}\label{rem_classification}
	By \Autoref{rem:ostriksresult}, all possible categorifications of the Grothendieck rings $K_n$ for $n>2$ and $K'_{n}$ for any $n$ can be found as the even part of the categorifications of Verlinde fusion rings, i.e., there are no new categories. We would like this to hold also for $K_0$, $K_1$, and $K_2$.
	
	Indeed, by a case-by-case distinction, we verify this. We will however, see that the number of possible categorifications does not follow the formula the Lemma would suggest, hence one can not simply extend it. For $n=0,1,2$, respectively, a naive interpolation of the lemma would state that there are $2,4,4$ categorifications, respectively.
	
	For $n=0$, the basis of the fusion ring $K_0$ has one element. There is only one fusion category categorifying it, namely the trivial fusion category $\mathrm{Vec}$. 
	
	For $n=1$, the fusion ring $K_1$ is isomorphic to the Grothendieck ring of $K(\mathrm{Vec}(\mathbb{Z}/2\mathbb{Z}))$, which has two categorifications. Both of them, however, can be found as the even part of some $\mathcal{C}_3^l$. The deviation from \Autoref{rem:ostriksresult} comes from the fact that while for $l\in\{1,3,5,7\}$ the categories $\mathcal{C}_3^l$ are pairwise non-equivalent, their even parts only give two non-equivalent categories. 
	
	Only the case $n=2$ remains. Following \cite[Chapter 5]{etingof2004classification}, any categorification of the fusion ring $K_2$ is a so-called \emph{group-theoretical} fusion category with classifications also found in the same source. Here we also have only two possible categorifications. However, both of them are also found as fusion subcategories in $\mathcal{C}_5^l$.
	
	We conclude therefore that the categories $\mathrm{Ad}(\mathcal{C}_n^l)$ are exactly those that categorify the even part of Verlinde fusion rings. We can distinguish different categories also by calculating only one value. For a fixed fusion ring, the dimension of $X_1$, i.e., $[3]$, tells us exactly which categorification we are looking at. Hence, it suffices to calculate these.
\end{remark}

\subsection{The asymptotic Hecke category for dihedral groups}
The two $J$-cells of size $1$, $c_0=\{1\}$ and $c_p=\{w_0\}$, have only one possible fusion categorification, as there is only one fusion category with one object, the category of finite dimensional vector spaces $\mathrm{Vec}$. The asymptotic Hecke category therefore is this trivial category, we can label its simple object by $B_1$ or $B_{w_0}$ depending on which cell we focus on. 

For the middle cell one can do diagrammatic calculations to see that the associators coincide exactly with the ones from $A_n$-fusion rule categories. We do not want to introduce all the notation for diagrammatic Soergel bimodules. Alternatively we use a result of Ben Elias that gives a connection of the diagrammatic category in the dihedral case to that of Temperley--Lieb algebras. 

\begin{theorem}\label{rem:elias}
	Let $W=I_2(p)$ be the dihedral group of order $2p$, let $c$ be the middle cell and $h\subset c$ a diagonal $H$-cell. Then the asymptotic Hecke category $\mathcal{H}_W^h$ is equivalent to $\mathrm{Ad}(\mathcal{C}_{n}^1)$ for $n=p-1$.
\end{theorem}
\begin{proof}
	First note that the $J$-ring of $h$ is the same as the Grothendieck ring of \Autoref{rem:ostriksresult}. Here $j_s$ corresponds to $X_0$, $j_{sts}$ to $X_1$, and so on. Hence, if $p$ is even we get the Grothendieck ring $K_{\frac{p-2}{2}}$; if $p$ is odd we get $K'_{\frac{p-3}{2}}$.

	As explained at the end of \Autoref{rem_classification}, it is therefore enough to calculate just one structure constant of the asymptotic category. By \cite{elias2015}, the (two-colored) Temperley--Lieb category embeds as the degree 0 morphisms into the category of Soergel bimodules of a dihedral group. By \cite[Theorem 2.15]{Elias_2017} we even have a degree-zero equivalence. 
	
	This equivalence between the diagrammatic dihedral and Temperley--Lieb calculations shows that the morphism spaces in the asymptotic Hecke category, and therefore the associators and dimensions, are described in exactly the same way. The computation of the associator data, which was described for the Temperley--Lieb category in \cite{edie-morrison-17}, works in the same way in the quotient of the diagrammatic Hecke category. As seen in \Autoref{lem:classification_ak}, the important structure constant is the choice of root of unity for $[2]$, i.e., the evaluation of the bubble in the Temperley--Lieb category. Under the equivalence of \cite{elias2015}, the bubble corresponds to the dimension of $B_{st}$ inside the asymptotic category, i.e., the composition $B_s\to (B_{st} \otimes B_{ts})(-1) \xrightarrow{\eta} (B_{st} \otimes B_{ts})(+1) \to B_s$. Since by construction $\partial_s(\rho)=1$, the composition evaluates to $\partial_s(\alpha_t)$. 
	
	The convention of \cite{IntroSoergel} there was to choose $[2]=-a_{s,t}=2\cos(\frac{\pi}{m_{s,t}}),$ hence we are in the case $\mathcal{C}_n^1$.
\end{proof}


The associator data computation works in the same way in the uncolored Temperley--Lieb category as done by Kauffman and Lins, see \cite{kauff1994}. However, in $\mathcal{H}_W^c$ we have two simple objects of any given length as a word can start in two different reflections. Note, that what is called \emph{$6j$-symbol}, i.e. the numerical data encoding the associator data, see \cite[Example 4.6.3]{EGNO}, therefore appears twice in the categorical data of $\mathcal{H}_W^c$, once for each color of the starting letter. Note further, that by \cite[Theorem 9.22]{IntroSoergel} the idempotent of $B_{w}$, with $w=st\ldots$ of length $k+1$, comes from the Jones--Wenzl idempotent of length $k$ inside the Temperley--Lieb category. Hence, the decrement of length by $1$. This is enough to describe the complete fusion structure of the asymptotic Hecke category associated to the middle cell in dihedral type.

\begin{corollary}
	Let $n+1=p\geq 3$ and consider the Coxeter group $W$ of type $I_2(p)$. Let $c$ be the subregular cell. The asymptotic Hecke category $\mathcal{H}_W^c$ associated to $c$ has the following fusion data.\begin{itemize}
		\item The objects are labelled by elements of the cell: $B_{w}$ for $w\in c$.
		\item The monoidal product is as in \Autoref{eq:mult_jring}, where $j_x$ denotes the equivalence class of $B_x$ in the Grothendieck ring.
		\item The associators are of the same form as for $A_n$-fusion rule categories, see \Autoref{ex:twistedZZ2ZZ}, where for an object $B_x$ we use the associators from $X_{\mathrm{len}(w)-1}$.
	\end{itemize}
\end{corollary}

To explain the last point: if $(B_x,B_y,B_z)$ has non-zero product, $B_w $ is a summand of $B_x \star B_y \star B_z$ and $B_u$ respectively $B_v$ are summands of $B_x \star B_y$ respectively $B_y \star B_z$, and we set $a,b,c,d,e,f$ to be the respective lengths minus $1$ of $u,v,w,x,y,z$, then the $6j$-symbol of $\mathcal{H}_W^c$ will be the same as that of the summands $X_a$ and $X_b$ for the summand of $X_c$ in $X_d \otimes X_e \otimes X_f$ inside the category $\mathcal{C}_n$.

\subsection{The center of $A_n$-fusion rule categories}\label{sub:center}

We can now investigate the center of the asymptotic Hecke category of the dihedral group by considering the categories $\mathrm{Ad}(\mathcal{C}_n)$.

First, we describe the Drinfeld center of $\mathcal{C}_n$. The main idea is to find a braiding on the category as this gives an equivalence to the center. The result will be independent on the chosen braiding on the category.

\begin{lemma}\label{lemma:center_of_braided} 
	Let $\mathcal{C}$ be a braided fusion category with invertible $S$-matrix. The center of $\mathcal{C}$ has the form \begin{equation}
		\mathcal{Z}(\mathcal{C}) \simeq \mathcal{C}\boxtimes \mathcal{C}^{rev},
	\end{equation} where $(-)^{rev}$ denotes the category $\mathcal{C}$ with reverse braiding, i.e. $c'_{X,Y}=c_{Y,X}^{-1}$ for $c_{X,Y}: X\otimes Y \to Y \otimes X$, the braiding on $\mathcal{C}$. 
	\begin{proof}
		This result is originally by Mueger \cite{mueger2003}, see also \cite[Propositions 8.6.1 and 8.20.12 ]{EGNO}. They show that the functors $\mathcal{C}\to\mathcal{Z}(\mathcal{C}),~X\mapsto (X,c_{-,X})$ and $\mathcal{C}^{rev}\to\mathcal{Z}(\mathcal{C}),~X\mapsto (X,c_{X,-}^{-1})$
		combine into an equivalence of braided tensor functors \begin{equation}
			\mathcal{C}\boxtimes\mathcal{C}^{rev}\to\mathcal{Z}(\mathcal{C}).
		\end{equation} 
		Note, that the center does not depend on the braiding chosen on $\mathcal{C}$ as long as the associated $S$-matrix is invertible. Hence, we can freely choose the braiding for computing the modular data of the center.
	\end{proof}
\end{lemma}

\begin{lemma}{\cite[Proposition 8.2.6]{FroehlichKerler1993}}\label{lem:classification_braiding}
	Let $1\leq l' \leq 4(n+1)$ be coprime to $n+1$ and $1\leq l \leq n$ such that $l \equiv \pm l'$ modulo $2(n+1)$. Then we can endow $\mathcal{C}_n^l$ with a braiding such that on the $X_2$ summand of $c_{X_1,X_1}: X_1\otimes X_1 \to X_1 \otimes X_1$ it is defined as $e^{\frac{2\pi il'}{4(n+1)}}$. This then gives a complete list of all braidings on a category with $A_{n}$-fusion rules.
\end{lemma}

Note, that one can compute the braidings of all other objects based on the one value given. Exact formulas are again found in \cite{edie-morrison-17}.

\begin{lemma}
	The $S$-matrix of the braided category $\mathcal{C}_{n}^l$ is 	\begin{equation}\label{eq:smatrix_typean}
		S_n=\left((-1)^{i+j}[(i+1)(j+1)]\right)_{i,j}.
	\end{equation} Here the indices $i,j$ are short for the objects $X_i, X_j$.
	\begin{proof}
		This is \cite[Equation 8.2.98]{FroehlichKerler1993} also computed in  \cite[Section 9.9]{kauff1994}.
	\end{proof}
\end{lemma}

\begin{example}\label{ex:center_c_3}
	One can check that the $S$-matrix of $\mathcal{C}_{n}^l$ is always invertible. Hence,
	\Autoref{lemma:center_of_braided} tells us directly that $\mathcal{Z}(\mathcal{C}_n)$ has $n^2$ simple objects. The object $X_i\boxtimes X_j$ maps to a simple object $X_i\otimes X_j$ in $\mathcal{Z}(\mathcal{C}_n)$ with a certain braiding coming from $X_i\in\mathcal{C}_n$ and $X_j\in\mathcal{C}_n^{rev}$.
	
	In the Deligne tensor product we get the $S$-matrix to be $S_n\otimes \bar{S_n}=S_n\otimes S_n$, i.e. the Kronecker product of the matrix with itself.  

	We can visualize this for $n=3$ in the following way. The $3^2=9$ simple objects of $\mathcal{Z}(\mathcal{C}_3)$ under the forgetful functor to $\mathcal{C}_3$ can be arranged  into a grid: \begin{equation}\label{eq:center_a_n}
		\begin{matrix}
			X_0\boxtimes X_0 & 	X_0\boxtimes X_1 & 	X_0\boxtimes X_2 \\ 
			X_1\boxtimes X_0 & 	X_1\boxtimes X_1 & 	X_1\boxtimes X_2 \\
			X_2\boxtimes X_0 & 	X_2\boxtimes X_1 & 	X_2\boxtimes X_2 \\
		\end{matrix} \qquad \rightsquigarrow \qquad 
		\begin{matrix}
			X_0 & X_1 & X_2 \\
			X_1 & X_0\oplus X_2 & X_1 \\
			X_2 & X_1 & X_0 \\
		\end{matrix}.
	\end{equation} The left side shows objects (without braiding) in $\mathcal{C}_3\boxtimes \mathcal{C}_3^{rev}$, the right side depicts the corresponding object in $\mathcal{Z}(\mathcal{C}_3)$ under the forgetful functor. We see that $X_1$ occurs with four different braidings in $\mathcal{Z}(\mathcal{C}_3)$, while $X_0$ and $X_2$ only with two. Furthermore, there is a simple object in $\mathcal{Z}(\mathcal{C}_3)$, namely $X_0\oplus X_2$ after applying the forgetful functor, which is obviously not simple in $\mathcal{C}_3$. Note also that all objects in $\mathcal{Z}(\mathcal{C}_3)$ are self-dual as they are self-dual in $\mathcal{C}_3$ and hence also in the Deligne tensor product. 
	
	The $S$-matrix of $\mathcal{C}_3$ is of the form \begin{equation}\label{eq:S3}
		S_3 = \begin{pmatrix}
			[1] & -[2] & [3] \\
			-[2] & [4] & -[6] \\
			[3] & -[6] & [9] \\
		\end{pmatrix} = \begin{pmatrix}
			1 & -\sqrt{2} & 1 \\
			-\sqrt{2} & 0 & \sqrt{2} \\
			1 & \sqrt{2} & 1\\ 
		\end{pmatrix}.
	\end{equation}
\end{example}

\begin{remark}[Lusztig's $S$-matrix]
	The dihedral fusion datum by Lusztig, \cite[Section 3.10]{lus94-exotic}, is of the following form: for $p\geq 3$ we consider the pairs $(i,j)$ with $0<i<j<i+j<p$ or $0=i<j<\frac{p}{2}$, as well as two special tuples $(0,\frac{p}{2})$ and $(0,\frac{p}{2})'$ if $p$ is even. We then define a pairing via \begin{equation}\label{eq:pairing_lusztig}
		\langle (i,j), (k,l) \rangle \coloneqq \frac{\xi^{il+jk}+\xi^{-il-jk}-\xi^{ik+jl}-\xi^{-ik-jl}}{p}
	\end{equation} on non-special tuples. Here $\xi$ is $e^{2\pi\cdot i / p}$, a primitive $p$-th root of unity. This expression looks similar to an expression in quantum numbers, the connection has been described in \cite[Section 3.4]{Lacabanne_2020}. 
	
	We set $n\coloneqq p-1$, then the tuples $(i,j)$ correspond to the object $X_{j-i-1}\boxtimes X_{j+i-1}$ in $\mathcal{C}_n\boxtimes \mathcal{C}_n^{rev}$. In the adjoint part both special elements will then correspond to two different subobjects of $X_{\frac{n-1}{2}}\boxtimes X_{\frac{n-1}{2}}$, see \Autoref{ex:i25}.
	
	For any tuple of pairs $((i,j),(k,l))$ the $S$-matrix value of the corresponding entry of $(X_{j-i-1}\boxtimes X_{j+i-1},X_{k-l-1}\boxtimes X_{k+l-1})$ is then \begin{equation}
		(-1)^{j+k-i-l-2}[(j-i)(k-l)][(j+i)(k+l)].
	\end{equation} The quantum part of this expression then gives \begin{align}
		&\frac{q^{(j-i)(k-l)}-q^{-(j-i)(k-l)}}{q-q^{-1}}\frac{q^{(j+i)(k+l)}-q^{-(j+i)(k+l)}}{q-q^{-1}} \\
		&=\frac{(q^{kj-ik-lj+il}-q^{ik-kj+lj-il})(q^{jk+jl+ik+il}-q^{-jk-jl-ik-il})}{(q-q^{-1})^2} \\
		&=\frac{q^{2kj+2il}-q^{-2ik-2jl}-q^{2ik+2lj}+q^{-2il-2jk}}{(q-q^{-1})^2},
	\end{align} where $q$ is a $2(n+1)$-th root of unity such that $q^2=\xi$. Indeed, this gives the result of the pairing by Lusztig modulo a factor of the form $\frac{(q-q^{-1})^2}{p}$, which is exactly the square root of the categorical dimension as in \Autoref{rem:normalize}.
\end{remark}

\subsection{The center of $\mathrm{Ad}(\mathcal{C}_n)$}

Here we describe the Drinfeld center of $\mathrm{Ad}(\mathcal{C}_n)$ as calculated by \cite{ediemichell2018brauerpicard}. We put it together with \cite{Lacabanne_2020} to compute its $S$-matrix and see how the normalized $S$-matrix is the same matrix Lusztig computed under in \cite[Section 3]{lus94-exotic} under an involution, i.e. a permutation on the columns. There is a case distinction depending on the parity of $n$. 

\subsubsection{The case of $n = 2m$ even}\label{subsubs_evenn}

It was noted in \cite[Lemma 3.1]{ediemichell2018brauerpicard} that the braiding of $\mathcal{C}_n$ restricted to the adjoint part $\mathrm{Ad}(\mathcal{C}_n)$ is still modular, i.e. the corresponding $S$-matrix is still invertible. In this case we can use \Autoref{lemma:center_of_braided} again.

\begin{lemma}[{\cite[Lemma 3.1]{ediemichell2018brauerpicard}}]\label{lemma:center_a_n_even}
	For $n = 2m$ even, we have as braided categories
	\begin{equation}
		\mathcal{Z}(\mathrm{Ad}(\mathcal{C}_{2m})) \simeq \mathrm{Ad}(\mathcal{C}_{2m}) \boxtimes \mathrm{Ad}(\mathcal{C}_{2m}^{rev}).
	\end{equation}
\end{lemma} 

\begin{example}\label{example_fibonacci}
	For $p = 5$ (the dihedral group $I_2(5)$), we have $n=4,~m=2$. The adjoint part of $\mathcal{C}_4$ is the Fibonacci category $\mathcal{F}$. We have two simple objects $(X_0,X_2)$ with monoidal product $X_2\otimes X_2\simeq X_0\oplus X_2$ and trivial associators except for triple $(X_2,X_2,X_2)$. Here, the associator on the second summand of $X_0 \oplus X_2^2 \simeq (X_2\otimes X_2) \otimes X_2 \mapsto X_2 \otimes (X_2\otimes X_2) \simeq X_0 \oplus X_2^2$ is the map 
	\begin{equation}
		X_2^2\to X_2^2,~\begin{pmatrix}
			\frac{[1]}{[3]} & -\frac{[2]^2}{[4]} \\ -\frac{[4]}{[2]^2[3]} & \frac{[6]}{[3][4]}
		\end{pmatrix} = \begin{pmatrix}
			\varphi^{-1} & -1-\varphi \\ -\varphi^{-3} & -\varphi^{-1}
		\end{pmatrix},
	\end{equation} 
	where $\varphi = \frac{1+\sqrt{5}}{2}$ and $[i]$ are the quantum numbers with $[2]=\varphi$.
	
	Furthermore, the $S$-matrix is the restriction of the $S$-matrix of $\mathcal{C}_4$, $S_4$, to the rows and columns corresponding to the adjoint part:
	\begin{equation}
		S_\mathcal{F}=\begin{pmatrix}
			[1] & [3] \\
			[3] & [9] \\
		\end{pmatrix}=\begin{pmatrix}
			1 & \varphi \\ \varphi & -1
		\end{pmatrix}.
	\end{equation} 
	Note that $[9]=-[1]=-1$ and $[3]=\varphi$. This is invertible, as expected by \Autoref{lemma:center_a_n_even}.
	
	Applying the forgetful functor we can visualize the objects of the center $\mathcal{Z}(\mathrm{Ad}(\mathcal{C}_{4})) = \mathcal{Z}(\mathcal{F})$ as the black objects in the matrix $X_i\otimes X_j$, see \Autoref{eq:center_a_n}
	\begin{equation}\label{eq:center_fibonnaci}
		\begin{matrix}
			X_0 & \cdot & X_2 & \cdot \\
			\cdot & {\color{gray}X_0\oplus X_2} & \cdot & {\color{gray}X_2}  \\
			X_2 & \cdot & X_0\oplus X_2 & \cdot \\
			\cdot   & {\color{gray}X_2} & \cdot &  {\color{gray}X_0} \\
		\end{matrix}
	\end{equation}
	This gives the $S$-matrix of $\mathcal{Z}(\mathcal{F})$ to be \begin{equation}
		S_\mathcal{F}\otimes \overline{S_\mathcal{F}} = 	\varphi\begin{pmatrix}
			\varphi^{-1} & 1 & 1 & \varphi \\
			1 & -\varphi^{-1} & \varphi &- 1 \\
			1 & \varphi & -\varphi^{-1} & -1 \\
			\varphi & -1 & -1 & \varphi^{-1}\\
		\end{pmatrix}.
	\end{equation} Here the ordering of objects is following the columns of \Autoref{eq:center_fibonnaci}, i.e. $X_0,X_2,X_2$, and then $X_0\oplus X_2$. 
	
	This matrix corresponds to Lusztig's result in \cite[Section 3.10]{lus94-exotic} under reordering and normalizing by the square root of the dimension of $\mathcal{Z}(\mathcal{C}_4)$ and applying an involution as seen in \cite[Remark before Proposition 3.1]{Lacabanne_2020}. To be more precise, we have $\dim(\mathcal{Z}(\mathcal{C}_4))=\dim(\mathcal{C}_4)^2$. To normalize the $S$-matrix we hence divide by $\dim(\mathcal{C}_4)=\dim(X_0)^2+\dim(X_2)^2=1^2+\varphi^2=\frac{5+\sqrt{5}}{2}=\sqrt{5}\varphi$. Under the ordering $(X_0,X_0\oplus X_2,X_2,X_2)$ we then get \begin{equation}
		\frac{1}{\sqrt{5}}\begin{pmatrix}
			\varphi^{-1} & \varphi & 1 & 1 \\
			\varphi & \varphi^{-1} & -1 &- 1 \\
			1 & -1 & -\varphi^{-1} & \varphi \\
			1 & -1 & \varphi & -\varphi^{-1}\\
		\end{pmatrix}.
	\end{equation} The final twist comes from the involution $(-)^\flat$, which sends $(i,j)\mapsto (i,p-j)$ if $i\geq 0$ and is trivial otherwise, see \cite[Section 3.1]{lus94-exotic}. This interchanges both copies of $X_2$ (the ones coming from the pairs $(1,2)$ and $(1,3)$) and leaves the other two elements invariant. Under the involution we therefore get exactly the matrix of \cite[Section 3.10]{lus94-exotic}.
\end{example}

\begin{remark}\label{rem:normalize}
	The calculations from \autoref{example_fibonacci} work generally for any $n=2m$ even, see the calculations in \cite[Section 3.4]{Lacabanne_2020}. We have $n^2$ simple objects in $\mathcal{Z}(\mathcal{C}_n)$ and hence $m^2$ in $\mathcal{Z}(\mathrm{Ad}(\mathcal{C}_n))$. The values of the normalized $S$-matrix coincide with the calculations done in \cite{Lacabanne_2020}. 
\end{remark}

\subsubsection{The case of odd $n$}\label{subsub_oddn}

Now, we consider the category $\mathrm{Ad}(\mathcal{C}_{2n+1})$. Here the restriction of the $S$-matrix is not invertible anymore, for example in \Autoref{eq:S3} the odd rows and columns give $\begin{pmatrix}
	1 & 1 \\ 1 & 1 
\end{pmatrix}.$ 

Therefore, one cannot use \Autoref{lemma:center_of_braided} directly. There is an alternative way described in \cite[Section 3]{ediemichell2018brauerpicard}.

\begin{construction}
	For a fusion category $\mathcal{C}$ and fusion subcategory $\mathcal{D}$ we denote the relative center $\mathcal{Z}_{\mathcal{D}}(\mathcal{C})$ as in \cite[Section 2.2]{gelaki2009centers}. Let $G$ be a finite group.
	
	If $\mathcal{C}$ is a $G$-graded fusion category, the trivial component $\mathcal{D}\coloneqq \mathcal{C}_0\subseteq\mathcal{C}$ is a fusion subcategory. By \cite[Theorem 3.5]{gelaki2009centers}, we have an equivalence \begin{equation}
		\mathcal{Z}_{\mathcal{D}}(\mathcal{C})^G \simeq \mathcal{Z}(\mathcal{C}).
	\end{equation}With this we can recover $\mathcal{Z}(\mathcal{D})$ from $\mathcal{Z}(\mathcal{C})$. The simple objects in $\mathcal{Z}(\mathcal{C})$ restricting to direct sums of the monoidal unit in $\mathcal{C}$ under the forgetful functor $\mathcal{Z}(\mathcal{C})\to\mathcal{C}$ form a subcategory $\mathcal{E}\simeq \mathrm{Rep}(G)\subseteq \mathcal{Z}(\mathcal{C})$. We get an equivalence \begin{equation}
		(\mathcal{E}')_G \simeq \mathcal{Z}(\mathcal{D}),
	\end{equation} where $(-)_G$ stands for the de-equivariantization.
	\begin{proof}
		This is \cite[Construction 3.2]{ediemichell2018brauerpicard} using \cite[Section 2 and 3 and Corollary 3.7]{gelaki2009centers}
	\end{proof}
\end{construction} 

\begin{example}\label{ex:i25}
	We continue with \Autoref{ex:center_c_3}, i.e. $n=3$. Here the categories $\mathcal{C}_n$ are $G\coloneqq \ZZ/2\ZZ$-graded, and $\mathcal{D}\coloneqq \mathrm{Ad}(\mathcal{C}_n)$ is the even or adjoint part of the category.
	
	The subcategory $\mathcal{E}\simeq \mathrm{Rep}(\ZZ/2\ZZ)$ has two simple objects, which we can understand as $X_0 \boxtimes X_0$ (the monoidal unit) and $X_{n-1} \boxtimes X_{n-1}$. Under the forgetful functor both are $X_0$, however with different braidings.
	The braidings on the second object are trivial on $X_0$ and $X_2$ but non-trivial on $X_1$, specifically: \begin{equation}
		c_{X_{0},X_1}: X_0\otimes X_1 \simeq X_1\to X_1 \simeq X_1 \otimes X_0
	\end{equation}
	is $-\id_{X_1}$. From this we can compute the centralizer $\mathcal{E}'$. Since the half-braiding of $X_{n-1} \boxtimes X_{n-1}$ on $X_1$ is non-trivial, no copy of $X_1$ can lie in $\mathcal{E}'$ as their half-braidings on the objects in $\mathcal{E}$ are trivial. All other objects however lie in the centralizer, i.e.
	all black objects in \begin{equation}
		\begin{matrix}
			X_0 & {\color{gray}X_1} & X_2 \\
			{\color{gray}X_1} & X_0\oplus X_2 & {\color{gray}X_1} \\
			X_2 & {\color{gray}X_1} & X_0 \\
		\end{matrix}.
	\end{equation}
	
	Under the de-equivariantization both copies of $X_0$ and $X_2$ in the corners will be isomorphic, the object $X_0\oplus X_2$ however decomposes into two simple objects $X_0$ and $X_2$ not isomorphic to the others. In total, we then get $4$ simple objects in the center $\mathcal{Z}(\mathrm{Ad}(\mathcal{C}_n))$.
	
	The restriction of the $S$-matrix of $\mathcal{Z}(\mathcal{C}_n)$ to the objects $X_0,X_2$ and $X_0\oplus X_2$ has the form \begin{equation}
		\begin{pmatrix}
			1 & 1 & 2 \\ 1 & 1 & -2 \\ 2 & -2 & 0\\
		\end{pmatrix},
	\end{equation} the $S$-matrix of $\mathcal{Z}(\mathrm{Ad}(\mathcal{C}_n))$ is of the form \begin{equation}\label{smat_dieder}
		\begin{pmatrix}
			1 & 1 & 1 & 1 \\ 1 & 1 & -1 & -1 \\ 1 & -1 & 1 & -1 \\ 1 & -1 & -1 & 1 \\
		\end{pmatrix} 
	\end{equation}and we see that the sum of the third and fourth rows and columns is the same as before.
	
\end{example}

\begin{remark}\label{rem_lacabanne}
	The calculations regarding the centralizer of the $\mathrm{Rep}(\ZZ/2\ZZ)$ subcategory have been done by Lusztig in \cite[Section 3.8]{lus94-exotic} and, in more detail including the $S$-matrix in \cite[Section 3.4]{Lacabanne_2020}. We have $n^2$ objects in the center of $\mathcal{C}_n$. For $n$ odd, i.e. $n=2m+1$ we have $2m^2+2m+1$ objects in $\mathcal{E}$. Under the de-equivariantization we get isomorphisms from the objects $X_i\boxtimes X_j$ to $X_{m-i}\boxtimes X_{m-j}$. The object $X_{m}\boxtimes X_m$ decomposes into a direct sum of two simple objects in $\mathcal{Z}(\mathrm{Ad}(\mathcal{C}_n))$, hence we are left with $m^2+m+2$ simple objects, as it was conjectured. 
\end{remark}

\begin{theorem}\label{thm:dihedral}
	\Autoref{lusztig_conjecture} holds for type $W=I_2(p)$.
	\begin{proof}
		This follows from \Autoref{rem:elias} and the computations of this Section.
		To reiterate: We have seen that the adjoint part of $A_n$-fusion rule (here $n=p-1$) categories are equivalent to the asymptotic Hecke category $\mathcal{H}^h$ for $h$ a diagonal $H$-cell in type $I_2(p)$, see \Autoref{rem:elias}. Specifically if $h=\{s,sts,ststs,\ldots\}$, the identification of objects of $\mathrm{Ad}(\mathcal{C}_n)$ to $\mathcal{H}^h$ is via \begin{equation}
			X_0 \leftrightarrow B_s,~X_2\leftrightarrow B_{sts},~X_4\leftrightarrow B_{ststs},\ldots .
		\end{equation}
		The equality of the $S$-matrices and the Fourier matrix constructed by Lusztig was seen for the even case in \Autoref{rem:normalize} and for the odd case it follows from the calculations referred to in \Autoref{rem_lacabanne}. In \cite[Section 3.4]{Lacabanne_2020} it was shown that the $S$-matrix of the center of $\mathrm{Ad}(\mathcal{C}_{2n+1})$ is indeed the Fourier matrix by Lusztig.
	\end{proof}
\end{theorem}

Explicit values of associators or $S$-matrices can be calculated using the Software described in \cite{MaurerThiel-Center,Maurer-TensorCategories}.

\section{The types $H_3$ and $H_4$}\label{sec:h34}

We give an overview of the possible $S$-matrices occurring for Drinfeld centers of asymptotic Hecke algebras corresponding to $J$-cells in non-crystallo\-graphic finite Coxeter groups. The missing two types $H_3$ and $H_4$ are discussed and complemented by the works of the previous sections.

\subsection{Type $H_3$ and $H_4$}
All cells, their $a$-values and asymptotic Hecke algebras of the Coxeter groups $H_3$ and $H_4$ have already been computed. See for example \cite{alvis08} for data on $H_4$. It turns out that the diagonal $H$-cells occurring are nearly always rather small, having only one or two elements. In these cases we have mostly only one possible categorification, hence the corresponding $S$-matrices are easy to list. In a couple of cases the associator is not known, and more calculations are needed. However, combinatorial results by Brou\'e and Malle, see \cite[Section 7]{brouemalle93} tell us which categorification should be the right one assuming \Autoref{lusztig_conjecture} is true.

These observations have been made in \cite[Section 8]{Mackaay2023} by Mackaay, Mazorchuk, Miemitz, Tubbenhauer and Zhang. The asymptotic Hecke category associated to an $H$-cell is denoted there by $\mathcal{A}_\mathcal{H}$ and called \emph{asymptotic bicategory}. The construction can be found in \cite[Section 3.2]{Mackaay2023}. We reuse their results on asymptotic Hecke categories in types $H_3$ and $H_4$ and augment their observations by possible $S$-matrices.

Only in type $H_4$ there is one cell with a considerably bigger $H$-cell. The $J$-cell of $a$-value $6$ contains diagonal $H$-cells of sizes $14$, $18$ and $24$. A description of the asymptotic Hecke category associated to it is unknown, there is however a combinatorial result \cite{malle94} about the $S$-matrix of its center, which contains $74$ simple objects, assuming \Autoref{lusztig_conjecture} is true in this case. Note, that this means that the centers of the three different categorifications $\mathcal{H}_W^h$ of sizes $14,18$ and $24$ all need to be equivalent to a category with $74$ simple objects.

\begin{example}\label{ex:datah3h4}
In type $H_3$ we have seven $J$-cells with data

	\begin{equation}
	\begin{tabular}{c||c|c|c|c|c|c|c}
		$c$ label (artificial) & $1$ & $2$ & $3$ & $4$ & $5$ & $6$ & $7$ \\
		\hline \hline
		$\abs{c}$ & $1$ & $18$ & $25$ & $32$ & $25$ & $18$ & $1$\\
		\hline
		$a$-value & $0$ & $1$ & $2$ & $3$ & $5$ & $6$ & $15$ \\
		\hline
		size of diagonal $H$-cell & $1$ & $2$ & $1$ & $2$ & $1$ & $2$ & $1$ \\
		\hline
		asymptotic Hecke category & (A) & (B) & (A) & (C) & (A) & (D) & (A) \\
	\end{tabular}
	\end{equation} 

in type $H_4$ there are $13$ $J$-cells with data

{\tiny
	\begin{equation}
	\begin{tabular}{c||c|c|c|c|c|c|c|c|c|c|c|c|c}
		$c$ label & $1$ & $2$ & $3$ & $4$ & $5$ & $6$ & $7$  & $8$ & $9$ & $10$ & $11$  & $12$ & $13$\\
		\hline \hline
		$\abs{c}$ & $1$ & $32$ &  $162$ & $512$ & $625$ & $1296$ & $9144$ & $1296$ & $625$ & $512$ & $162$ & $32$ & $1$ \\
		\hline
		$a$-value & $0$ & $1$ &  $2$ & $3$ & $4$ & $5$ & $6$ & $15$ & $16$ & $18$ & $22$ & $31$ & $60$ \\
		\hline
		$\abs{h}$ & $1$ & $2$ & $2$ & $2$ & $1$ & $1$ & $14,18,24$ & $1$ & $1$ & $2$ & $2$ & $2$ & $1$ \\
		\hline
		$\mathcal{H}_W^h$ & (A) & (B) & (D) & (C) & (A) & (A) & (E) & (A) & (A) & (C) & (D) & (D) & (A)\\
	\end{tabular} 
	\end{equation}
} with the following cases for the asymptotic Hecke category:

\begin{enumerate}\label{enum_h4cases}
 	\item\label{itema} There is only one element in $h$, hence any categorification has only one simple object, we therefore have $\mathcal{H}_W^h\simeq \mathrm{Vec}$. 
 
 	\item\label{itemb} The fusion ring structure is that of the Fibonacci category, in the same way as in \Autoref{example_fibonacci}. We have a categorification through $\mathrm{Ad}(\mathcal{C}_4)$.
 
	\item\label{itemc} The fusion ring structure is that of $K(\mathrm{Vec}(\ZZ/2\ZZ))$. This has two different categorifications, we have either $\mathcal{H}_W^h\simeq \mathrm{Vec}(\ZZ/2\ZZ)$, or $\mathcal{H}_W^h\simeq\mathrm{Vec}^\omega(\ZZ/2\ZZ)$ for the non-trivial associator as in type $E_8$. 
	
	\item\label{itemd} Here we have the same fusion ring as in case (B), however it is not clear which root of unity is in the categorification. Either  we have $[2]=\frac{1+\sqrt{5}}{2}$ as in case (B) or $[2]=\frac{1-\sqrt{5}}{2}$. 
	
	\item This is the only case where no categories with the respective Grothen\-dieck ring are known. If \Autoref{lusztig_conjecture} holds we will have that for any such category $\mathcal{H}_W^h$ the center has $74$ simple objects.
\end{enumerate}
\end{example}

\subsection{The $S$-matrices of centers of asymptotic Hecke categories of exotic cells}

 We complete the overview of $S$-matrices of the centers of asymptotic Hecke categories associated to two-sided Kazhdan--Lusztig cells which was started in \Autoref{cor_smatrix_weyl} and extended in \Autoref{cor_e78}, using the results for type $H$ above and for the dihedral group as in \Autoref{thm:dihedral}.

\begin{theorem}\label{thm:overview}
	Let $c$ be a two-sided cell in a Coxeter group of type $I_2(p),H_3,\allowbreak{}H_4$. The $S$-matrix of the center of the asymptotic Hecke category is of the following type:
	\begin{enumerate}
		\item If $c$ contains a diagonal $H$-cell of size one, such as the cells $\{1\}$, $\{w_0\}$ in type $I_2(p)$, or the ones of case (A) in type $H_3,H_4$, the asymptotic Hecke category is $\mathrm{Vec}$ and the $S$-matrix is \begin{equation}
			S= \begin{pmatrix}
				1
			\end{pmatrix}.
		\end{equation}
	
		\item For the cases (B) in types $H_3,H_4$ we get the asymptotic Hecke category to be $\mathcal{H}_W^h\simeq \mathrm{Ad}(\mathcal{C}_4)$. The normalized $S$-matrix of its center is \begin{equation}
		S_\mathcal{F}=\frac{1}{\sqrt{5}}\begin{pmatrix}
				\varphi^{-1} & \varphi & 1 & 1 \\
				\varphi & \varphi^{-1} & -1 &- 1 \\
				1 & -1 & -\varphi^{-1} & \varphi \\
				1 & -1 & \varphi & -\varphi^{-1}\\
			\end{pmatrix},
		\end{equation} where $\varphi=\frac{1+\sqrt{5}}{2}$, see \Autoref{example_fibonacci}. 
	
		\item If $c$ is middle $J$-cell in type $W=I_2(p)$ the asymptotic Hecke category is $\mathcal{H}_W^c\simeq \mathrm{Ad}(\mathcal{C}_{n})$, for $n\coloneqq p-1\geq 2$. For $n$ even the $S$-matrix is \begin{equation}\label{eq:smatdiheral}
			S(\mathcal{Z}(\mathrm{Ad}(\mathcal{C}_n))) = \left([(2i-1)(2j-1)] \right)_{1\leq i,j \leq n}^{\otimes 2},
		\end{equation} for $[2]=q+q^{-1}$ and $q$ a $2(n+1)$-th root of unity. For the normalized $S$-matrix we divide by $\frac{(q-q^{-1})^2}{n+1}$. For $n=4$ this is exactly the result of the previous item. For $n$ odd we have seen that in $\mathcal{Z}(\mathcal{H}_{I_2(p)}^h)\simeq \mathcal{Z}(\mathrm{Ad}(\mathcal{C}_n))$ one simple object of $\mathrm{Ad}(\mathcal{C}_n)\boxtimes \mathrm{Ad}(\mathcal{C}_n)$ splits in the center, and the $S$-matrix therefore includes the matrix \eqref{eq:smatdiheral} as well as two new rows and columns, whose entries can be computed by the pairing of Lusztig, see \eqref{smat_dieder} in \Autoref{ex:i25} and \eqref{eq:pairing_lusztig}. For $I_2(4)$ this gives for example $S(\mathcal{Z}(\mathrm{Vec}(\ZZ/2\ZZ)))$ and for $I_2(6)$ we get $S(\mathcal{Z}(\mathrm{Vec}(S_3))$ see \Autoref{cor_smatrix_weyl}.
		
		\item\label{case:h} In the cases (D) in type $H_3,H_4$ we had two possible categorifications. The $S$-matrix will either be $S_{\mathcal{F}}$, or a variant of it, where we replace $[2]=\varphi$ by $\varphi^{-1}=\frac{1-\sqrt{5}}{2}$. We call this modified $S$-matrix $S_{\mathcal{F}'}$. In \cite[Equation 7.3]{brouemalle93} all Fourier matrices associated to these cases in type $H_3,H_4$ are of the form $S_{\mathcal{F}}$. If \Autoref{lusztig_conjecture} is true we therefore expect to never see the second option.

		\item In the case (C) in type $H_3,H_4$ we again had two possible categorifications, namely the category of $\ZZ/2\ZZ$-graded vector spaces, either with trivial or non-trivial twist. The possible normalized $S$-matrices are \begin{equation}
			S_C\in \Bigg\{\frac{1}{2}\begin{pmatrix}
				1 & 1 & 1 & 1 \\
				1 & 1 & -1 & -1 \\
				1 & -1 & 1 & -1 \\
				1 & -1 & -1 & 1 \\
			\end{pmatrix},\frac{1}{2}\begin{pmatrix}
				1 & 1 & 1 & 1 \\
				1 & 1 & -1 & -1 \\
				1 & -1 & -1 & 1 \\
				1 & -1 & 1 & -1 \\
			\end{pmatrix} \Bigg\}.
		\end{equation} However, if \Autoref{lusztig_conjecture} holds the Fourier matrix as computed by \cite{brouemalle93} is the second option, i.e. the same as in the exceptional cases of type $E_7$ and $E_8$, see \Autoref{cor_e78}.
	
		\item And finally we expect the $S$-matrix for the cell of $a$-value $6$ in type $H_4$ to be the Fourier matrix computed by \cite{malle94}. 
	\end{enumerate} 
\end{theorem}

\section{Examples in infinite Coxeter groups} \label{sec:infinite_cases}
So far we have only considered finite Coxeter groups. In these cases it is clear that all two-sided cells themselves are also finite. However, one might also investigate finite two-sided cells lying in infinite Coxeter groups.

There are conjectured results on the structure of cells in infinite Coxeter groups, see \cite{belolipetsky2014cells,belolipetsky2014kazhdanlusztig}. However, as far as the authors are aware, there are no classification results on finite two-sided or $H$-cells.

We will focus this subsection on two known classifications and extend the description of $S$-matrices of asymptotic Hecke categories to all finite two-sided cells of $a$-value lower than or equal to $2$. 

The cell of $a$-value $0$ is always finite as it only contains the unit. The asymptotic Hecke category in this case is the category of finite-dimensional vector spaces $\mathrm{Vec}$, the $S$-matrix of its center is $(1)$. 

In general, let $W$ always denote a Coxeter group with generating set $S$. We write $W_i\coloneqq \{x\in W\mid a(x)=i\}$ for the subsets of elements of a given $a$-value. We say that $W$ is \emph{$a(i)$-finite} if $W_i$ is finite.

\subsection{The case of $a(1)$-finite Coxeter groups}

 Hart gave a characterization of $a(1)$-finite Coxeter groups in \cite{hart2017elements}. The set $W_1$ has always an easy description.

\begin{remark}\label{rem_cella1}
	In any irreducible Coxeter group $W$, the set $W_1$ forms a two-sided cell, often referred to as the \emph{subregular cell}. This is shown, for example, in \cite[Chapter 12]{Bonnafe18}. The elements of $W_1$ are precisely those elements of $W$ that have a unique reduced expression (also established in \cite[Chapter 12]{Bonnafe18}). Furthermore, the left and right cells within $W_1$ are determined by the left and right descent sets of the elements, i.e. for $x, y \in W_1$, we have $x \sim_L y$ if and only if the unique reduced expressions of $x$ and $y$ end in the same reflection. A complete description of the corresponding $J$-ring is given in \cite{Xu2019}.
\end{remark}

\begin{theorem}\label{lem_a(1)-finite}
	Let $(W,S)$ be a Coxeter system. Then $W_1$ is finite if and only if $S$ is finite and each connected component of $S$ is a tree with no infinite bonds and has at most one edge label greater than three.
	\begin{proof}
		This is \cite[Corollary 2.6]{hart2017elements}.
	\end{proof}
\end{theorem}

We give a short overview of the combinatorial proofs of \cite[Section 2]{hart2017elements} as this helps us to describe the asymptotic Hecke category. Using the following lemma one can reduce the theorem to the study of irreducible Coxeter groups:

\begin{lemma}\label{lem:addition_a}
	Let $(W,M)$ be a Coxeter system and let $K,L\subset M$ be disjoint subsets of $M$ not sharing a connected component inside the Dynkin diagram, i.e. $m_{s,t}=2$ for $s\in K,t\in L$. 
	We denote the Coxeter groups generated by $K$ and $L$ by $U$ and $V$, then $U\times V \subset W$ is a subgroup of $W$. For two-sided cells $c_1\subset U$ and $c_2\subset V$ of $a$-value $i$ and $j$ the Cartesian product $c\coloneqq c_1\times c_2\subset W$ lies inside a two-sided cell of $a$-value $i+j$ inside $W$. The asymptotic Hecke algebra $J_{c}$ is also isomorphic to $J_{c_1}\times J_{c_2}$.
	\begin{proof}
		This follows quickly from the observation that for $x\in c_1$ and $y\in c_2$ the Kazhdan--Lusztig basis elements commute, i.e. we have $b_xb_y=b_yb_x$ and therefore $b_{(x,y)}=b_xb_y$. The cell and $a$-value computations now work independently in both summands.
	\end{proof} 
\end{lemma}

\begin{remark}
	Let now $(W,S)$ be an irreducible Coxeter system. Then the main steps in the proof of \Autoref{lem_a(1)-finite} are:
	\begin{itemize}
		\item Note that an expression $(s_1,s_2,\ldots,s_n)$ is never unique if there is an $i$ with $m_{i,i+1} = 2$. Therefore, all unique expressions must be paths inside the Dynkin diagram.
		\item If $S$ is infinite the elements of $\{s \mid s \in S\}$ all lie in the subregular cell, so we can assume $S$ to be finite. Also, any cycle inside the Dynkin diagram gives rise to infinitely many elements of $a$-value $1$. The same happens for $m_{s,t}=\infty$.
		
		
		\item Finally let $(s,t)$ and $(u,v)$ be tuples in the same connected component with $m_{s,t},m_{u,v}>3$ and take a path $p$ connecting both tuples. Without loss of generality, we have $p=(t=r_0,r_1,\ldots,r_n=u)$ and neither $s$ nor $v$ occur inside $p$. Let $p^{-1}$ denote the reverse path. Then the composition $(p,v,p^{-1},s)$ represents a reduced expression and any power of it does too, hence $W_1$ is infinite. We need $m_{s,t},m_{u,v}$ to be greater than $3$ as otherwise the composition contains a subsequence $(u,v,u)$ or $(t,s,t)$ not unique.
	\end{itemize}	
\end{remark}

	The observation on paths can also be used to enumerate the elements of $W_1$, describe the left and right cell structure as well as give the multiplication table of the $J$-ring.

\begin{corollary}\label{cor:enumeration_a1-finite}
	Let $W$ be an $a(1)$-finite irreducible Coxeter group and let $c\coloneqq W_1$. Let $m$ be the value of the biggest relation occurring in the Dynkin diagram. For any tuple $(r,s)$ of reflections in $W$ there is a unique $H$-cell $h_{r,s}$ where all words start in $r$ and end in $s$. The size is $\lfloor \frac{m}{2} \rfloor$ if the shortest path connecting $r$ and $s$ includes the edge $m$ and $\lfloor \frac{m-1}{2} \rfloor$ if not.
	\begin{proof}
		This follows from the proof of \Autoref{lem_a(1)-finite} by counting the number of paths corresponding to reduced expressions. We need the characterization of left and right cells inside $W_1$ by starting and ending letter as seen in \Autoref{rem_cella1}. An enumeration of $W_1$ can also be found in \cite[Theorem 2.5]{hart2017elements}.
	\end{proof}
\end{corollary}

\begin{theorem}
	Let $W$ be an $a(1)$-finite Coxeter group and let $c\subseteq W_1$ be a two-sided cell. Then one can choose an $H$-cell $h\subset c$ such that the asymptotic Hecke category $\mathcal{H}_W^h$ is equivalent to $\mathrm{Ad}(\mathcal{C}_n)$ for some $n$, i.e. the center and the $S$-matrix are the same as in the dihedral case as seen in \Autoref{thm:overview}.
	\begin{proof}
		We assume that $W$ is irreducible. Let $s,t$ be generators of $W$ such that $m_{s,t}$ is maximal (i.e. we take the unique tuple $(s,t)$ such that $m_{s,t}$ is greater than $3$ if it exists). We now choose the $H$-cell starting and ending in $s$, $h\coloneqq h_{s,s}$. This cell is then the same as the $H$-cell of the subgroup generated only by $s$ and $t$, a dihedral group of order $2m_{s,t}$. All computations of $\mathcal{H}_W^h$ therefore reduce to the finite dihedral case.
	\end{proof}
\end{theorem}

\begin{example}
	One such Coxeter group has appeared in \cite[Figure 1]{belolipetsky2014kazhdanlusztig}. The Coxeter group is of type $W_{237}$  with generators $\langle r,s,t \mid r^2=s^2=t^2=(rs)^3=(st)^7=(rt)^2=1\rangle$. Following \Autoref{cor:enumeration_a1-finite} we can enumerate all elements of the subregular cell by looking for paths corresponding to reduced expressions. We order these elements by starting and ending letter, i.e. we partition them into left and right cells:
	\begin{table}[!ht]
		\centering
		\begin{tabular}{c|c|c}
			$\{r,rstsr, rststsr\}$ & $\{sr,stsr, ststsr\}$ & $\{tsr,tstsr,tststsr\}$ \\
			\hline
			$\{rs,rsts,rststs\}$ & $\{s,sts,ststs\}$ & $\{ts,tsts,tststs\}$ \\
			\hline
			$\{rst,rstst,rststst\}$ & $\{st,stst,ststst\}$ & $\{t,tst,tstst\}$ \\
		\end{tabular}
	\end{table}
	On the diagonal $H$-cell coming from the dihedral subgroup of type $I_2(7)$ the multiplication on the asymptotic Hecke algebra can be read of directly. We have for example $j_{sts}^2=j_s+j_{sts}+j_{ststs}$. Similarly, one can work out the complete multifusion ring structure and get for example $j_{sr}j_{rststsr}=j_{ststsr}$. The center of the asymptotic Hecke category has $14$ simple objects.
\end{example}

\subsection{The case of $a(2)$-finite Coxeter groups}

Recent results by Green and Xu classified all irreducible Coxeter groups which are $a(2)$-finite. Coxeter groups where the Dynkin diagram contains a cycle have either none or infinitely many elements of $a$-value $2$. For all other cases they further always described one $H$-cell lying in $W_2$. We list their results and show that the $S$-matrix of the asymptotic Hecke category is the same as in the dihedral case of \Autoref{thm:overview}, by choosing an appropriate $H$-cell in which the asymptotic Hecke algebra is isomorphic to the Grothendieck ring of $\mathrm{Ad}(\mathcal{C}_n)$. 

\begin{proposition}{\cite[Theorem 3.31 and Theorem 4.17]{green2020classification} and \cite[Proposition 4.1]{green2020classification}} \label{prop:a2-finite}
	An irreducible Coxeter group $W$ with elements of $a$-value $2$ is $a(2)$-finite if and only if it is of one of the following types:\begin{equation}
		A_n,~B_n,~\tilde{C}_n,~E_{q,r},~F_n,~H_n,~I_n,
	\end{equation} 
	where 
	\begin{equation}
		\tilde{C}_{n-1}=\dynkin[extended,Coxeter,labels={1,2,n-1,n}, scale=1.8] C{o...oo},
	\end{equation}	
	for $n\geq 3$,
	 \begin{equation}
	 	E_{q,r}=\dynkin[%
	 	labels={-q,v,-(q-1),0,(r-1),r},
	 	scale=1.8] E{ooo...o...oo},
	 \end{equation}	for $r\geq q\geq 1$,
 	\begin{equation}
 		F_n=\dynkin[Coxeter,
 		labels={1,2,3,n},
 		scale=1.8] F{ooo...o},
 	\end{equation} for $n\geq 4$, 
 	\begin{equation}
 	H_n=\dynkin[Coxeter,
 	labels={1,2,3,n},
 	scale=1.8] H{ooo...o},
 \end{equation} for $n\geq 3$.
 
	In the case $E_{q,r}$ where $r=q=1$ (i.e. $D_{4}$) the set $W_2$ consists of three two-sided cells, if $r>q=1$ (i.e. type $D_{r+3}$) we have two cells in $E_{q,r}$. In all other cases $W_2$ itself is a two-sided cell. 
	
	One representative of an $H$-cell is given by the following:
	
	\begin{itemize}
		\item Type $A_n$: $h=\{13\}$
		\item Type $B_3$: $h=\{13\}$
		\item Type $B_n$ for $n>3$: $h=\{24,2124\}$
		\item Type $\tilde{C}_{n-1}$ where $n\geq 5$: $h=\{24,2124,2z,212z\}$, where $z=45\ldots\allowbreak(n-1)n(n-1)\ldots 54$
		\item Type $E_{q,r}$ where $r\geq q\geq 2$: $h=\{1v\}$
		\item Type $F_4$: $h=\{24\}$
		\item Type $F_n$, where $n>4$: $h=\{24,243524\}$
		\item Type $H_3$: $h=\{13\}$
		\item Type $H_n$, where $n>3$: $h=\{24,2124\}$
	\end{itemize}
\end{proposition}

\begin{theorem}
	Let $W$ be an irreducible $a(2)$-finite and infinite Coxeter group. The center of the asymptotic Hecke category associated to $c=W_2$ the two-sided cell of $a$-value $2$ is equivalent to one of the following cases:
	
	\begin{equation}
		\mathrm{Vec}, ~\mathcal{Z}(\mathcal{F}),~\mathcal{Z}(\mathcal{F}'), ~\mathcal{Z}(\mathcal{F})\boxtimes \mathcal{Z}(\mathcal{F}),
	\end{equation} for $\mathcal{F}$ the Fibonacci category $\mathrm{Rep}(\mathfrak{so}(3)_3)$ as in \Autoref{example_fibonacci} and $\mathcal{F}'$ the Fibonacci category with the second choice of associator as seen in \Autoref{thm:overview}. 
	
	The possible $S$-matrices are:
	
	\begin{equation}
		(1),~S_\mathcal{F}\otimes S_\mathcal{F},~S_{\mathcal{F}'}\otimes S_{\mathcal{F}'},~(S_\mathcal{F}\otimes S_\mathcal{F})\otimes (S_\mathcal{F}\otimes S_\mathcal{F}).
	\end{equation} 
	
	\begin{proof}
		From the classification in \Autoref{prop:a2-finite} we are only concerned with the infinite cases $\tilde{C}_n,E_{q,r},F_m,H_l$ for $q\geq 2$ and $r+q\geq 7$, $m>4,l>5$. If the $H$-cell given has size $1$ then the only possible categorification of the asymptotic Hecke algebra is $\mathrm{Vec}$.
		
		Therefore, we only need to check the other three remaining cases. In all of them we find that the $H$-cell lies in a finite parabolic subgroup in which it is also an $H$-cell of $a$-value $2$. In this parabolic subgroup $h$ can be written as $h_1\times h_2$, where both are of $a$-value $1$ inside the respective Coxeter group, see \Autoref{lem:addition_a}. One can therefore deduce the structure of the asymptotic Hecke algebra from finite cases.
		
		\begin{itemize}
			\item Case $F_m$ for $m>4$: the $H$-cell $h=\{24,243524\}$ lies in the parabolic subgroup $B_4$ where we identify the generators $i$ of $B_4$ by $i+1$ inside $F_4$. The asymptotic Hecke algebra in this case is the Fibonacci ring, where $j_{24}$ is the identity and $j_{243524}^2=j_{24}+j_{243524}$. Therefore, \Autoref{rem:ostriksresult} holds and we categorify over an $A_n$-fusion rule category as in the dihedral case.
			
			\item Case $H_n$ for $n>3$: the $H$-cell $h=\{24,2124\}$ lies in the finite parabolic subgroup $I_2(5)\times A_1$, where we identify the generators $1$ and $2$ with those of $I_2(5)$ and $4$ with the one of $A_1$. This case therefore also reduces to the observations of \Autoref{thm:overview}.
			
			\item Case $\tilde{C}_n$: the $H$-cell $h=\{24,2124,2z,212z\}$ lies in the parabolic subgroup $I_2(4)\times B_{n-3}$. It is the product of their unique cells of $a$-value $1$, namely $\{2,212\}$ in $I_2(4)$ and $\{4,z\}$ in $B_4$. Both asymptotic Hecke categories associated to them are again the Fibonacci category. The asymptotic Hecke category associated to $h$ is therefore the product $K(\mathcal{F})\times K(\mathcal{F})$. It could be that there are different categories categorifying this ring, however by a classification result on small rank fusion categories \cite[Theorem 1.1]{Liu_2022} only the Deligne tensor product of the Fibonacci category with itself lies over the given fusion ring. The asymptotic Hecke category $\mathcal{H}_W^h$ is therefore equivalent to $\mathcal{F}\boxtimes \mathcal{F}$, the center is $\mathcal{Z}(\mathcal{F}\boxtimes \mathcal{F})\simeq \mathcal{Z}(\mathcal{F})\boxtimes \mathcal{Z}(\mathcal{F})$ and the $S$-matrix is given by the Kronecker product as stated.
		\end{itemize} 
	\end{proof}
\end{theorem}

In these cases more diagrammatic calculations are needed to find which of the options is the right categorification.
\end{document}